\documentclass[review,12pt,authoryear]{elsarticle}

\usepackage[T1]{fontenc}
\usepackage[unicode]{hyperref}
\usepackage[margin=1in]{geometry}
\usepackage{amssymb}
\usepackage{amsmath}
\usepackage{float}
\usepackage{subcaption}
\usepackage{multirow}

\makeatletter
\def\ps@pprintTitle{%
	\let\@oddhead\@empty
	\let\@evenhead\@empty
	\def\@oddfoot{\footnotesize\itshape Preprint submitted to \@journal\hfill}%
	\let\@evenfoot\@oddfoot}
\makeatother

\begin{document}

\begin{frontmatter}
\journal{Journal}
\title{Improving Hydrofoil Efficiency Through Simulation-Driven Design Optimization}

\author[institution1,institution2]{Ante Sikirica}
\ead{ante.sikirica@uniri.hr}

\author[institution1]{Ivana Lu\v{c}in\texorpdfstring{\corref{corresponding}}{}}
\ead{ivana.lucin@riteh.uniri.hr}

\author[institution1]{Marta Alvir}
\ead{marta.alvir@riteh.uniri.hr}

\author[institution1]{Tin Matulja}
\ead{tin.matulja@riteh.uniri.hr}

\cortext[corresponding]{Corresponding author}

\affiliation[institution1]{organization={Faculty of Engineering, University of Rijeka},
	addressline={Vukovarska 58},
	city={Rijeka},
	postcode={51000 Rijeka},
	country={Croatia}}

\affiliation[institution2]{organization={Center for Advanced Computing and Modelling, University of Rijeka},
	addressline={Radmile Matejčić 2},
	city={Rijeka},
	postcode={51000 Rijeka},
	country={Croatia}}

\begin{abstract}
Foil windsurfing is a sport in which competitors fine-tune their setup to optimize performance under varying environmental conditions. However, such adjustments are typically guided by intuition. This paper presents a methodology for improving hydrofoil design using computational fluid dynamics. First, structured light scanning was used to obtain the geometry of a conventional hydrofoil. Subsequently, five derived designs have been explored, inspired by advancements in airplane wing design. Following this, a preliminary design was selected, and a design of experiments approach was implemented. Generated simulation data were used to train a machine learning surrogate, which is ultimately used in the multi-objective design optimization. While the surrogate-based approach yielded promising results, it was limited to exploring a restricted set of design parameters. To overcome this constraint, additional improvements were attained through adjoint optimization, demonstrating that comparable performance can be achieved with a distinctly different hydrofoil design. Finally, the conventional design was refined through adjoint optimization. Notably, adjoint optimization yielded a 15\% increase in efficiency. The results suggest that both proposed approaches can produce novel hydrofoil designs capable of outperforming the original configuration under specific conditions; however, due to its inherent flexibility, the adjoint approach with realistic constraints offers more freedom in the design exploration process.
\end{abstract}

\end{frontmatter}

\section{Introduction}

Foil windsurfing is a water sport developed from traditional windsurfing in which hydrofoils are mounted under surfboards to improve overall sailing performance. These hydrofoils typically consist of two wings connected by a fuselage and a mast that extends above the surface of the water. As the board moves through the water due to the wind acting on the sail, the hydrofoil generates lift, causing the board to rise above the water. This leads to a reduction in drag. Consequently, sailing speeds are higher, and maneuverability is improved. Foil windsurfing is a competitive sport with no strict equipment regulations. Competitors can, therefore, make modifications to their equipment based on their own experience or rely on innovative engineering methods and techniques to gain an advantage. In search of novel designs, biologically inspired concepts have attracted attention across various engineering fields. For hydrofoils, nature-inspired shapes such as the flippers of humpback whales have already been considered \citep{kant2020}. Likewise, advances in aerodynamics have shown that winglets and wingtip devices can improve performance \citep{guerrero2020}. Drawing on insights from these areas may also accelerate progress in hydrofoil design.

The numerical evaluation of different hydrofoil designs and configurations requires considering both aerodynamic and hydrodynamic effects. Capturing the dynamic behavior of the system additionally calls for unsteady simulations. However, computational fluid dynamics (CFD) simulations can be computationally expensive, especially in the early stages of design. Consequently, naval architects often rely on velocity prediction programs (VPPs), which are based on numerical, experimental, or empirical data. Although these programs were initially developed for traditional displacement vessels, they have since been adapted to accommodate a wider range of complex sailing conditions. Examples include the development of six-degree-of-freedom VPPs for racing yachts \citep{kerdraon2020}, performance prediction of hydrofoiling lake racers \citep{melis2022}, and for America's Cup foiling vessels \citep{patterson2022}. More recently, a model for windsurfing velocity prediction was developed \citep{hou2024} that considered the aerodynamic and hydrodynamic properties of the sail, board, and sailor.

CFD simulations have also been used to analyze the hydrodynamics of kitefoils \citep{ocana2017} and to optimize the design of a kiteboarding hydrofoil \citep{backas2016}. A study by Sakellariou et al. \citep{sakellariou2017} explored the use of CFD to improve the shape of a surfboard fin. The authors chose forty-two designs through Latin hypercube sampling (LHS) and constructed a response surface model (RSM) based on the CFD results. Ultimately, a genetic algorithm was employed to identify the optimal shape that achieved the desired performance target. Similarly, Wang et al. \citep{wang2018} employed a sampling strategy for CFD simulations to build a surrogate model, which was subsequently used to optimize the cross-section of an underwater glider. Sacher et al. \citep{sacher2018} combined a two-dimensional flow solver and a structural model to optimize a hydrofoil for the America's Cup. The paper by Pehlivan Solak et al. \citep{pehlivan2024} relied on CFD data to optimize the main wing of a kitefoil. A recent study by Alonso and Calderon-Sanchez \citep{alonso2025} presented a systematic analysis of the influence of key design features. The authors considered winglets, varying strut angles, bulbs, and proximity to the free surface and assessed how these features affect hydrodynamic performance. Nowruzi et al. \citep{nowruzi2017} explored the use of artificial neural networks (ANN) to predict the lift-to-drag ratio for two chosen NACA hydrofoils and different geometric characteristics, namely the aspect ratio (AR) and the taper ratio (TR), different Reynolds numbers, and angle of attack (AoA). The paper by Coppede et al. \citep{coppede2019} considered Gaussian processes as a surrogate in order to improve the design and consequently the resistance of the hull in calm water. Similarly, Kim et al. \citep{kim2024} employed deep neural networks as surrogate models in a design optimization process to predict and improve hull resistance. Recently, digital twins built around a baseline design have been employed to systematically explore new design variants. By combining simulations with experimental insights, they could be used to predict, test, and improve real-world performance, with computational resources as the primary additional investment \citep{liu2021, boyes2022}.

The extent of hydrofoil design improvement is fundamentally constrained by the manner in which the design space is defined. In other words, the efficiency and performance improvement achieved in any optimization process is linked to the employed parametrization. Conventional parametrization approaches, which typically utilize geometric characteristics such as chord length, span, tip angle, and aspect ratio, are interpretable but inherently restrict the flexibility of shapes that can be generated. Consequently, to explore unconventional designs, parametrization strategies must be more flexible and should not purely rely on fixed geometric descriptors.

Domain-based methods, such as free-form deformation, shift focus from geometric parameters to spatial deformation of control volumes. In this approach, a control volume encompassing the geometry or region of interest is first defined. The selected area is described using a lattice of control points whose movement deforms the entire enclosed mesh rather than modifying specific geometric features. When coupled with optimization algorithms, this approach enables large geometric deformations and can lead to innovative designs that might not emerge through conventional parametrization. However, while they offer greater geometric flexibility, they sacrifice the direct interpretability of conventional parameters, making it more challenging to establish clear relationships between design modifications and performance outcomes. Additionally, the increase in degrees of freedom can lead to non-physical geometries or geometries not suitable for manufacturing if not properly constrained. Nazemian and Ghadimi \citep{nazemian2021} employed a multi-objective (MO) optimization approach to improve the bow region of trimaran hulls, aiming to reduce total resistance in both calm water and waves. The hull shape was controlled using the arbitrary shape deformation technique with 120 design variables. The focus, however, was on the ship hull. A paper by Walker et al. \citep{walker2024} provides an overview of shape optimization approaches for hulls and hydrofoils, with particular emphasis on parametrization strategies. Notably, only one study on hydrofoil optimization was identified \citep{bonfiglio2018}, but even this work did not employ a free-form deformation approach, highlighting the limited attention this area has received.

The present paper proposes an integrated methodology that combines reverse engineering, numerical simulations, and design optimization, drawing inspiration from recent advancements in aerodynamics. The process begins with reverse engineering, where a structured light scanning technique is used to capture and reconstruct the geometry of an existing hydrofoil. By digitizing the actual hydrofoil, a precise digital representation is obtained, which serves as a reliable baseline for further analyses. Numerical simulations are subsequently employed to systematically evaluate the hydrodynamic behavior of the scanned hydrofoil under different operating conditions. The simulations are used to identify the performance characteristics of the current design and highlight areas where improvements might be possible. Building on these insights, the methodology explores alternative design variants that aim to improve specific performance characteristics. In the optimization stage, an automated process is introduced that combines a statistical sampling technique, CFD simulations, and an ANN-based surrogate model. The surrogate model approximates the relationship between design parameters and performance metrics, providing a mapping that can be used to efficiently explore the design space. Coupled with a multi-objective optimization approach, this enables rapid investigation of new geometries that balance key performance indicators. The optimized design obtained through this approach was subsequently refined using a free-form deformation technique. This fine-tuning process aims to further improve the hydrodynamic performance. In addition, the initial reconstructed design was subjected to free-form deformation to assess whether improvements could still be achieved for the nominally optimal solution. A key contribution of this study is the direct comparative analysis of two distinct optimization paradigms applied to the same hydrofoil system under identical conditions. While surrogate-based and adjoint-based optimization have each been applied separately in marine hydrodynamic problems, few studies have systematically compared their relative performance in a practical engineering context. This head-to-head comparison reveals important trade-offs.

The presented workflow aims to address two main objectives. First, it provides a systematic framework for understanding how specific design modifications influence hydrodynamic characteristics. Second, it establishes a practical approach for performance improvement by reducing reliance on physical prototyping, thereby accelerating the design process.

\section{Materials and methods}
\subsection{Baseline geometry}

The front and rear wings employed in this study were obtained using a structured light scanning method. Digital reconstruction was performed intentionally, as the wings had already been used and modified in situ; their geometry therefore differed from the original design. The mast and fuselage geometries were obtained from publicly available data.

Structured light scanning is a three-dimensional, contactless technique based on light pattern projection. Compared to coordinate measuring machines, structured light scanners are less accurate but offer greater flexibility and are non-invasive \citep{bernal2013}. The accuracy depends on factors such as projector quality, camera resolution, and the properties of the light pattern. Measurements were conducted using a GOM ATOS II Triple Scan Rev. 02 blue light scanner at the Center for Advanced Computing and Modelling, University of Rijeka, in an air-conditioned environment at a constant temperature of $22^{\circ}$C. The cameras were set to position 400. Schneider Cinegon 1.4/12-0946 lenses were used. The scanner was calibrated using a standard CP40 calibration panel, with results showing deviations of less than 0.1 pixels for all parameters, thus meeting the manufacturer's requirements. The point spacing during measurement ranged from $0.07 \leq d_p \leq 0.62$ mm. Reference points with a diameter of $d_{rp} = 3$ mm were used. To improve scanning in critical regions, particularly along edges, an anti-reflection film and spray were applied. During the measurement process, the hydrofoils were mounted on a turntable using a custom-made vice.

\begin{figure}[H]
    \centering
    \begin{subfigure}[b]{0.495\linewidth}
        \centering
        \includegraphics[width=\linewidth]{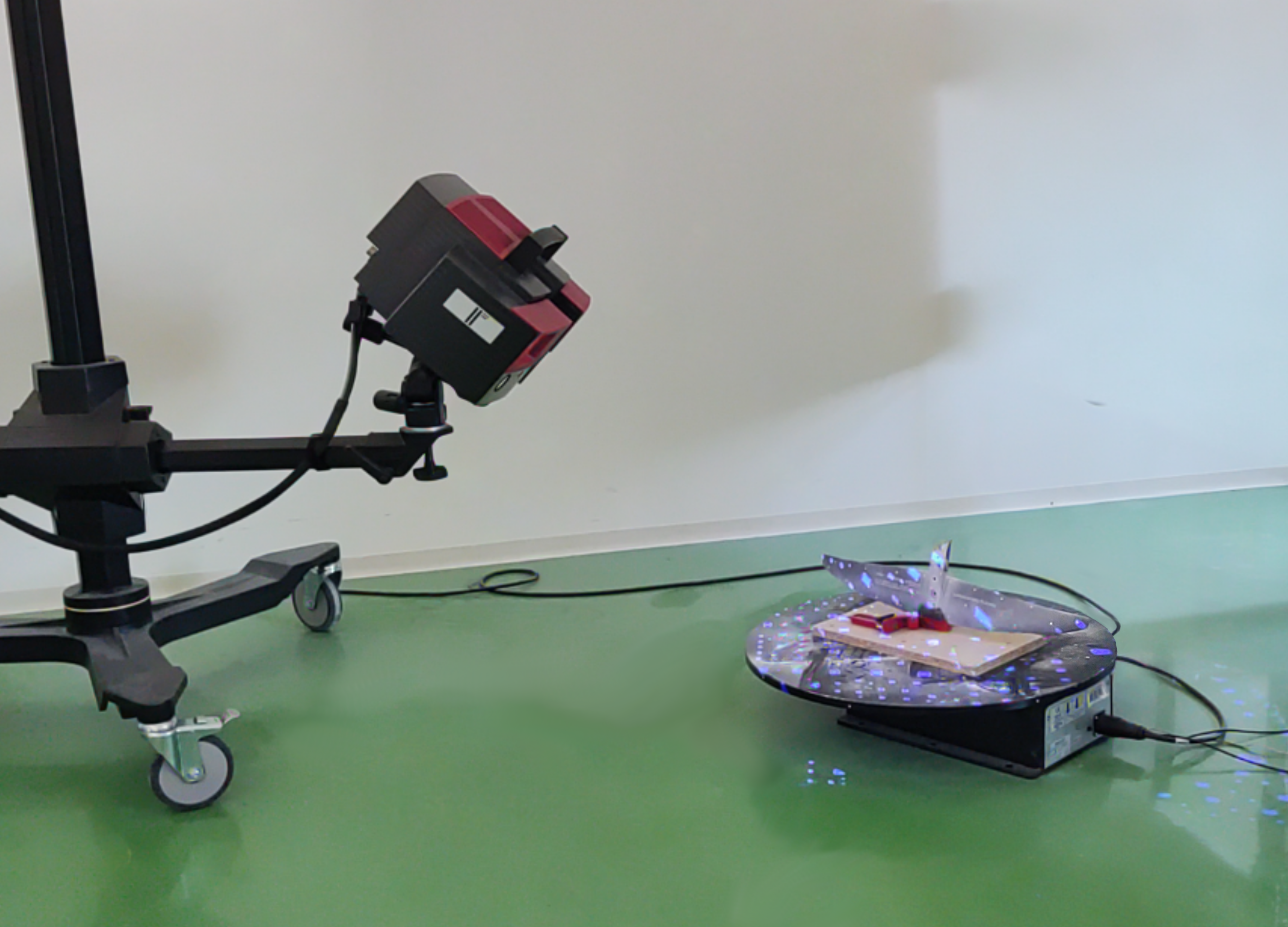}
        \caption{}
    \end{subfigure}
    \hfill
    \begin{subfigure}[b]{0.495\linewidth}
        \centering
        \begin{subfigure}[b]{\linewidth}
            \centering
            \includegraphics[trim={5cm 2.5cm 5cm 2.5cm}, clip, width=\linewidth]{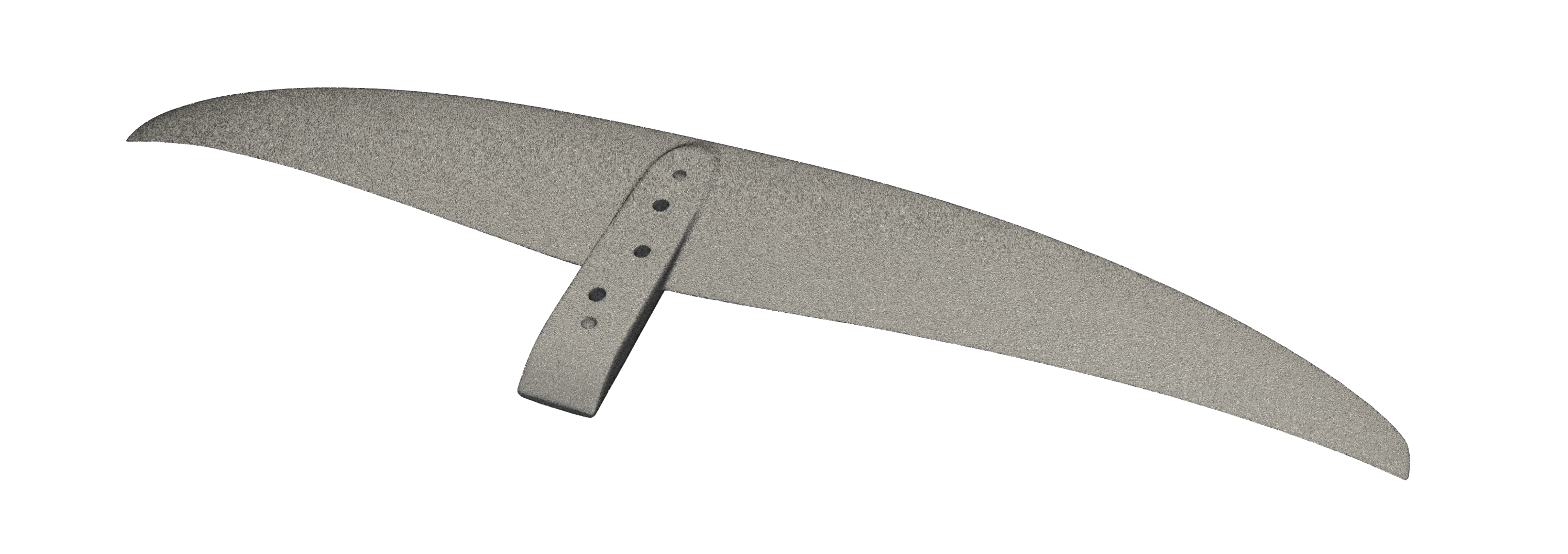}
            \caption{}
        \end{subfigure}
        \vfill
        \begin{subfigure}[b]{\linewidth}
            \centering
            \includegraphics[trim={5cm 5cm 5cm 5cm}, clip, width=\linewidth]{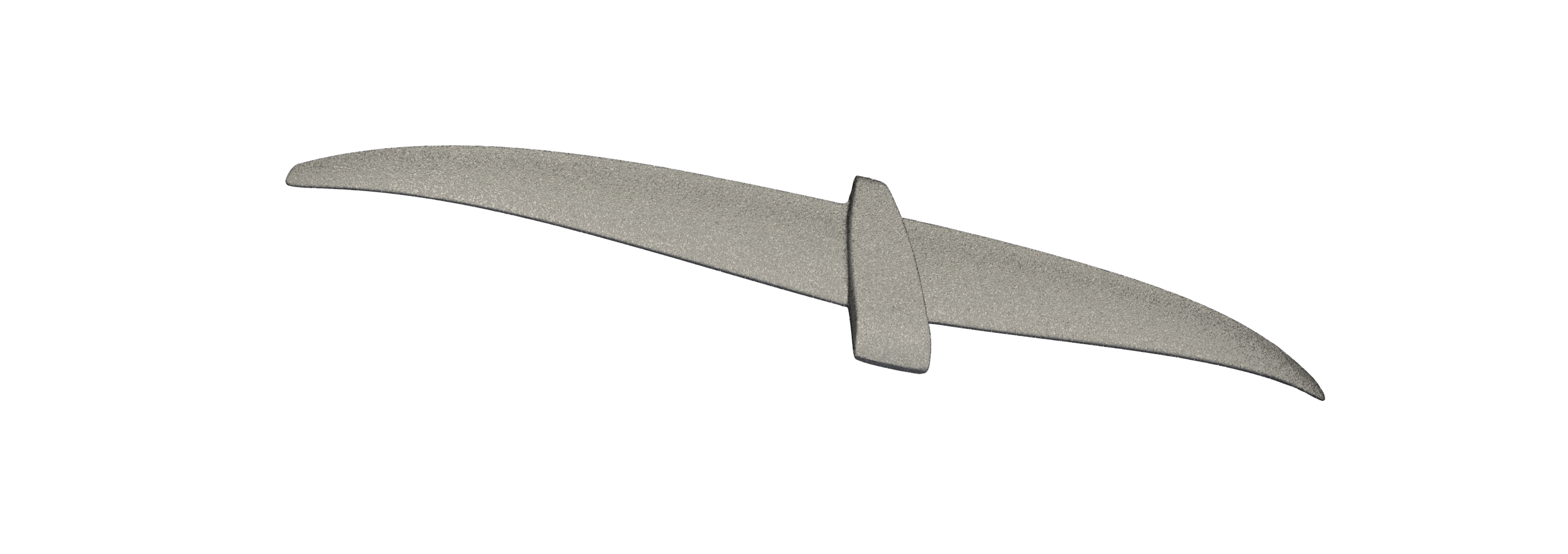}
            \caption{}
        \end{subfigure}
    \end{subfigure}
    \caption{Overview of the scanning process: (a) GOM ATOS II Triple Scan with turntable and front wing, (b) reconstructed front wing, and (c) reconstructed rear wing.}
    \label{fig:scanning_environment}
\end{figure}

The ATOS Professional V8 SR1 Suite was used for both data acquisition and initial geometry clean-up. The scanned data were exported as point clouds. Further processing and refinement of the point clouds were performed in Rhino 7, including adjustments to critical surface details. The final parameterized geometry, compliant with ISO 10303-21 standards, was created using the CATIA V5R20 CAD suite. The experimental setup and resulting geometries are shown in Figure \ref{fig:scanning_environment}.

\subsection{Preliminary configurations}
\label{sec:preliminary_configurations}

\begin{figure}[H]
    \centering
    \includegraphics[trim={20cm 1cm 20cm 1cm}, clip, width=0.9\linewidth]{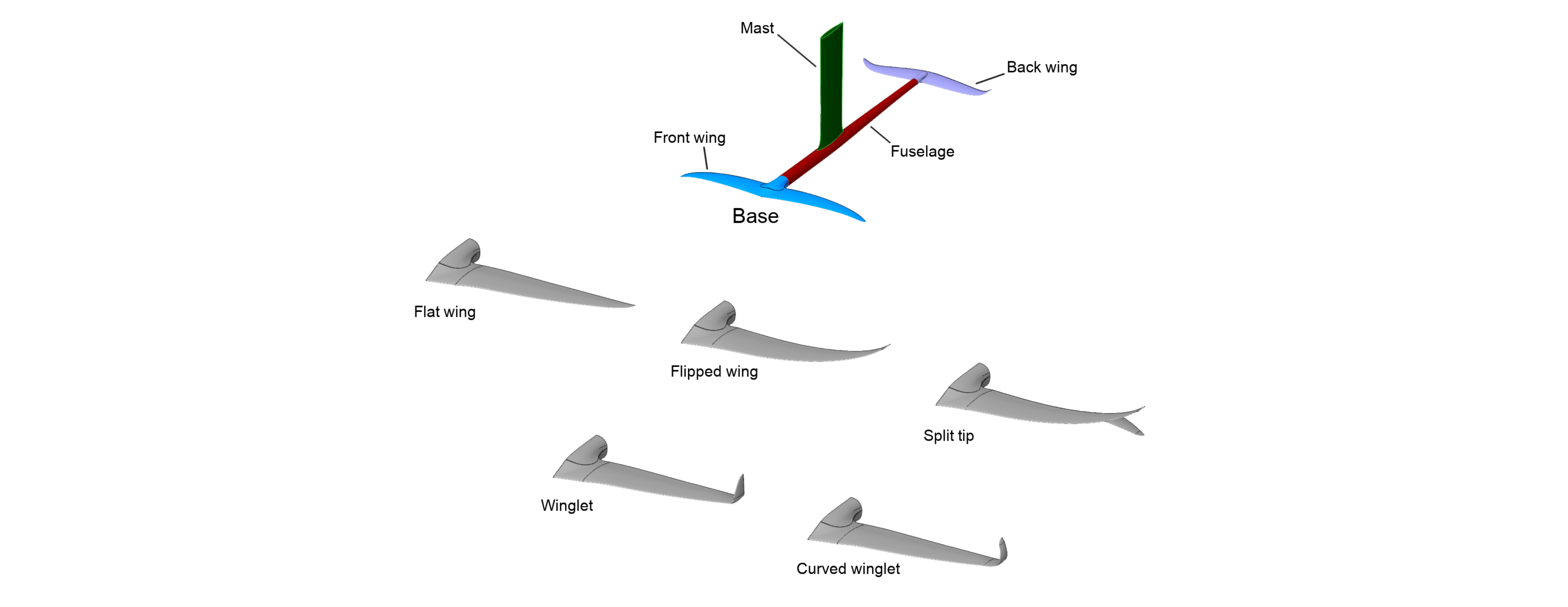}
    \caption{Overview of the considered front wing designs: base wing design and hydrofoil configuration (top), followed by modified front wing variants: flat wing, flipped wing, split tip, winglet, and curved winglet.}
    \label{fig:hydrofoil_designs}
\end{figure}

Given that the front wing exerts the greatest influence on overall fluid-flow characteristics, including drag and lift \citep{perali2024}, this study focused exclusively on front wing modifications, with no changes made to the rear wing, fuselage, or mast. In addition to the baseline design, five alternative front wing configurations were initially considered: a flat wing, a wing with upward-angled tip (analogous to aircraft wing dihedral), a split-tip winglet, a flat wing with a vertical winglet, and a flat wing with a curved winglet. Similar design concepts, including dihedral and winglet configurations, have been explored for aircraft wings \citep{guerrero2020}. Figure \ref{fig:hydrofoil_designs} shows the baseline design alongside all considered design variants.

Only the submerged portion of the hydrofoil was modeled. As the hydrofoil is assumed to be positioned 0.5 m below the water surface, the upper 0.35 m of the mast above the waterline was omitted. The total length of the hydrofoil is 0.932 m. For the baseline design, the front wing spans 0.58 m and the rear wing 0.40 m. The front wing has a maximum thickness of approximately 10 mm as it approaches the centerline, while the rear wing reaches about 7 mm at its thickest point.

It should be noted that all parameters of the considered design variants were chosen arbitrarily. The flat wing configuration is a simplified reference case, which isolates the effects of the baseline wing's anhedral geometry. The upward-angled flipped wing examines whether a dihedral approach might offer benefits in this context. The split-tip design aims to disrupt tip vortices, potentially reducing induced drag by dividing the primary vortex. Both the vertical and curved winglet configurations draw from conventional aircraft wing design, where winglets are widely used to mitigate induced drag by weakening wingtip vortices.

The geometries were designed to remain as close as possible to the base design in terms of overall size and projected frontal area. In this context, the wings can be viewed as deformed versions of the original geometry rather than entirely new designs. For instance, for the winglet and curved winglet configurations, despite the bent or curved tip geometries, the total frontal-projected area remains approximately equivalent to the baseline configuration. By preserving these key geometric parameters, the study isolates the influence of local shape modifications on lift and drag. This controlled strategy provides a basis for future parametric investigations.

\subsection{Numerical setup}

\begin{figure}[H]
    \centering
    \includegraphics[trim={15cm 0cm 20cm 0cm}, clip, width=\linewidth]{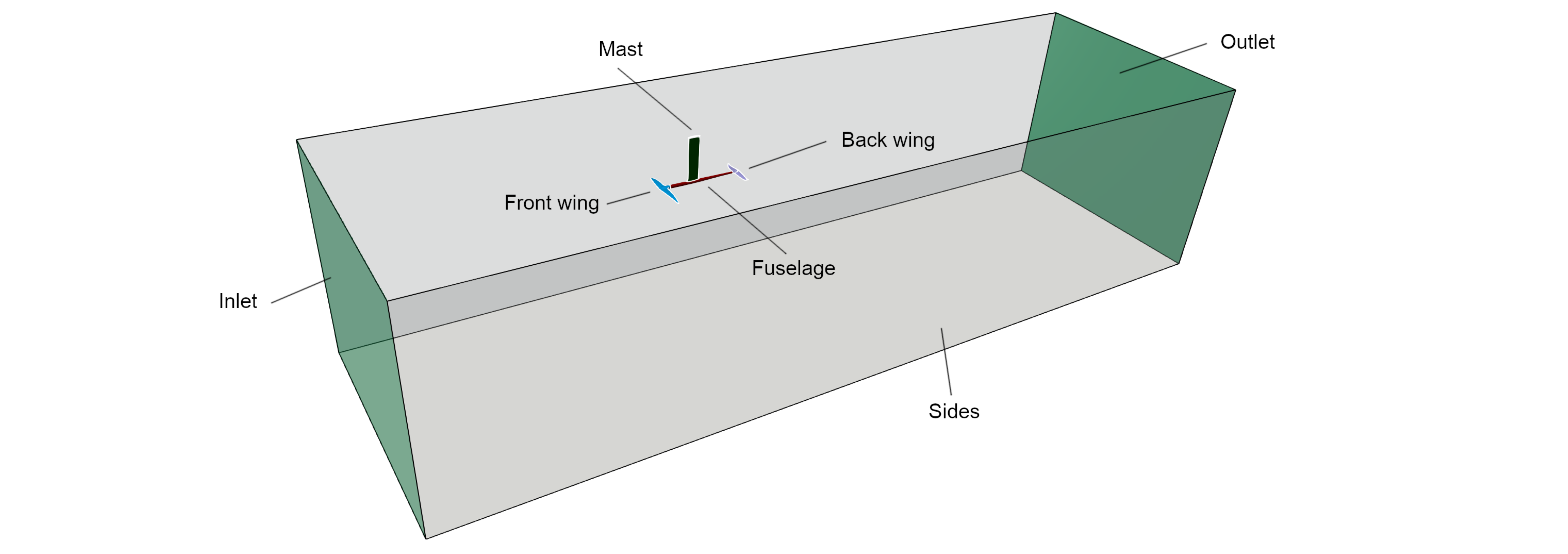}
    \caption{Computational domain and boundary patches.}
    \label{fig:domain}
\end{figure}

The computational domain for the numerical simulations measured 10 m in length, 3 m in width, and 2.5 m in height. The foil was placed 0.5 m below the top boundary and 3 m downstream of the inlet, ensuring sufficient upstream flow development and minimizing inlet effects. An illustration of the domain and boundary patches is shown in Figure \ref{fig:domain}.

The flow was assumed steady and incompressible. The Reynolds-Averaged Navier-Stokes (RANS) approach was used to model turbulence, with the equations for continuity and momentum defined as:
\begin{equation}
    \nabla \cdot \mathbf{u} = 0
\end{equation}
\begin{equation}
    (\mathbf{u} \cdot \nabla)\mathbf{u} = -\frac{1}{\rho}\nabla p + \nabla \cdot \left[(\nu + \nu_t)(\nabla \mathbf{u} + \nabla \mathbf{u}^T)\right] - \frac{2}{3}\nabla k
\end{equation}
where $\mathbf{u}$ is the mean velocity vector, $p$ the mean pressure, $\rho$ the fluid density, $\nu$ the kinematic viscosity, $\nu_t$ the eddy viscosity, and $k$ the turbulent kinetic energy. For turbulence closure, the SST $k$-$\omega$ model \citep{menter1994} was adopted, which introduces two additional transport equations for turbulent kinetic energy and the specific dissipation rate $\omega$:
\begin{equation}
    \nabla \cdot (\mathbf{u} k) = \nabla \cdot \left[(\nu + \sigma_k \nu_t)\nabla k\right] + P_k - \beta^* k\omega
\end{equation}
\begin{equation}
    \nabla \cdot (\mathbf{u} \omega) = \nabla \cdot \left[(\nu + \sigma_\omega \nu_t)\nabla \omega\right] + \gamma \frac{P_k}{\nu_t} - \beta \omega^2 + 2(1-F_1) \sigma_{\omega 2} \frac{1}{\omega} \nabla k \cdot \nabla \omega
\end{equation}
where $P_k$ is the production term, and the model constants follow standard SST $k$-$\omega$ formulation \citep{menter1994}.

The numerical modeling was performed using the OpenFOAM 10 toolbox \citep{of10}. The working fluid was water, with a kinematic viscosity of $\nu = 1 \cdot 10^{-6}$ m$^2$/s. A second-order accurate scheme, linearUpwind, was applied for the velocity term, with the upwind scheme used for the turbulence terms. Each simulation ran for up to $5 \times 10^4$ iterations or until the residuals for each variable dropped below $10^{-5}$. A summary of the prescribed boundary conditions is provided in Table \ref{tab:openfoam_bc}.

\begin{table}[H]
	\centering
	\footnotesize
	\caption{Employed boundary conditions.}
	\label{tab:openfoam_bc}
	\begin{tabular}{@{}lcccc@{}}
		\hline
		\textbf{Field} & \textbf{Inlet} & \textbf{Outlet} & \textbf{Sides} & \textbf{Hydrofoil patches} \\
		\hline
		$\mathbf{u}$ & fixedValue & pressureInletOutletVelocity & symmetry & fixedValue \\
		$p$ & zeroGradient & fixedValue & symmetry & zeroGradient \\
		$\nu_t$ & calculated & calculated & symmetry & nutkWallFunction \\
		$k$ & fixedValue & inletOutlet & symmetry & kqRWallFunction \\
		$\omega$ & fixedValue & inletOutlet & symmetry & omegaWallFunction \\
		\hline
	\end{tabular}
\end{table}

The influence of design changes was determined by measuring the forces and calculating the drag and lift coefficients. The drag coefficient ($C_d$) and lift coefficient ($C_l$) were chosen as objectives because they provide a direct measure of the hydrofoil's hydrodynamic performance. The primary purpose of a hydrofoil is to generate lift, raising the hull out of the water and thereby reducing drag; the lift coefficient quantifies this effect, while the drag coefficient measures the hydrodynamic cost associated with generating that lift \citep{faltinsen2005, anderson2011}. Together, their ratio defines the hydrodynamic efficiency, making them physically meaningful. Furthermore, both coefficients are nondimensional, enabling comparisons across different geometries, scales, and operating conditions. The drag coefficient is calculated as:
\begin{equation}
	C_d = \frac{F_d}{0.5 \, \rho A u^2}
\end{equation}
where $F_d$ is the drag force, $\rho$ is the fluid density, $A$ is the reference area, and $u$ is the inlet velocity. The lift coefficient is calculated analogously:
\begin{equation}
	C_l = \frac{F_l}{0.5 \, \rho A u^2}
\end{equation}
where $F_l$ is the lift force. The reference area used for all calculations was the planform-projected area of the front wing, rear wing, and fuselage at $0^\circ$ angle of attack. The inlet velocity was maintained at 20 m/s for all simulations. The designs were assessed at seven angles of attack: $-10^\circ$, $-5^\circ$, $0^\circ$, $5^\circ$, $10^\circ$, $15^\circ$, and $20^\circ$.

For the baseline design, a grid convergence study was carried out at $0^\circ$ AoA using three grid resolutions: coarse, medium, and fine. The maximum and minimum cell sizes for the fine grid were $l_{\mathrm{max}} = 0.036$ m and $l_{\mathrm{min}} = 0.00135$ m, respectively, with a first cell height of $l_{\mathrm{first}} = 0.00075$ m. The coarse and medium grids were generated by scaling the cell sizes with a factor of $f_s = 1.33$. The total number of cells in the coarse, medium, and fine grids was approximately $1.2 \cdot 10^6$, $2.1 \cdot 10^6$, and $4.0 \cdot 10^6$, respectively. The grid convergence index (GCI) \citep{roache1998} was calculated based on the drag coefficient, with results presented in Table \ref{tab:gci}. The values confirm that the solution is within the asymptotic range of convergence. Based on these findings, the fine grid was selected for all subsequent simulations. The average $y^+$ for the fine grid was approximately 120 on the front wing, with peak values of up to 250 on the mast and fuselage.

\begin{table}[H]
	\centering
	\footnotesize
	\caption{Results of the grid convergence study.}
	\label{tab:gci}
	\begin{tabular}{@{}lccccc@{}}
		\hline
		       & \boldmath$\varepsilon_{c,m}$ (\%) & \boldmath$\varepsilon_{m,f}$ (\%) & \boldmath$GCI_{c,m}$ (\%) & \boldmath$GCI_{m,f}$ (\%) & \boldmath$\frac{GCI_{m,f}}{r^p \cdot GCI_{c,m}}$ \\
		\hline
		$C_d$ & 15.20 & 5.04 & 8.78 & 2.91 & 0.95 \\
		\hline
	\end{tabular}
\end{table}

\subsection{Data-driven design}
\label{sec:data_driven_design}

Among the five preliminary designs (Section \ref{sec:preliminary_configurations}), the configuration with the most favorable lift-to-drag ratio was chosen for systematic parametric investigation. A design of experiments (DOE) study was conducted using full-factorial sampling to ensure adequate data coverage for surrogate model training.

All numerical simulations in this phase were performed with first-order discretization schemes to reduce computational cost. This simplification is justified, as the primary goal at this stage is to assess the relative performance trends among different design variants rather than to obtain high-fidelity quantitative results. In addition, minor differences in numerical accuracy are generally outweighed by the surrogate model's predictive accuracy and generalization capability. The final, optimal design was re-evaluated using second-order schemes to maintain consistency with the baseline methodology.

The data collected in the DOE phase were used to train a surrogate model. A multilayer perceptron (MLP) surrogate model was implemented using scikit-learn 1.5.1 \citep{scikit}. Hyperparameter tuning was conducted using Optuna 4.0 \citep{optuna}. The tuning process systematically explored a range of parameters, including the number of hidden layers (1 to 5), the number of neurons per layer (2 to 16), the activation function (ReLU or tanh), and the regularization parameter, $\alpha$, which was varied between $10^{-4}$ and $10^{-1}$ on a logarithmic scale. Additionally, the maximum number of epochs was included as a tunable parameter. The MLP was configured as a multi-output regressor to predict both the drag coefficient $C_d$ and lift coefficient $C_l$. Before training, all input features (DOE design variables) were scaled using a MinMaxScaler.

To ensure robustness and prevent overfitting, a nested validation strategy was adopted. The dataset was first partitioned into training (80\%) and holdout test (20\%) sets using stratified sampling to preserve representative distributions. During hyperparameter tuning, 5-fold cross-validation was used, and the weighted regression error metric $E_r$, a weighted sum of the component metrics was minimized:
\begin{equation}
    E_r = w_{R^2}\,(1 - E_{R^2}) + w_{MRE}\,E_{MRE} + w_{RRMSE}\,E_{RRMSE}
\end{equation}
where $E_{R^2}$ is the mean $R^2$ error, $E_{MRE}$ mean relative error, and $E_{RRMSE}$ mean relative root mean squared error, given by:
\begin{align}
    E_{R^2} &= \frac{1}{k} \sum_{j=1}^{k} \left( 1 - \frac{\sum_{i=1}^{n_j} (y_i^{(j)} - \hat{y}_i^{(j)})^2}{\sum_{i=1}^{n_j} (y_i^{(j)} - \bar{y}_{\vphantom{i}}^{(j)})^2} \right), \\
    E_{MRE} &= \frac{1}{k} \sum_{j=1}^{k} \left( \frac{1}{n_j} \sum_{i=1}^{n_j} \frac{|y_{i}^{(j)} - \hat{y}_{i}^{(j)}|}{|y_{i}^{(j)}|} \right), \\
    E_{RRMSE} &= \frac{1}{k} \sum_{j=1}^{k} \left( \frac{\sqrt{\frac{1}{n_j} \sum_{i=1}^{n_j} (y_{i}^{(j)} - \hat{y}_{i}^{(j)})^2}}{\bar{y}_{\vphantom{i}}^{(j)}} \right)
\end{align}
where $k$ is the number of folds, $n$ denotes the number of samples, $y_i$ the true values, $\hat{y}_i$ the predicted values, and $\bar{y}$ the mean of the true values.

The tuned MLP surrogate was subsequently coupled with the multi-objective optimization framework jMetalPy 1.6 \citep{jmetalpy}. The non-dominated sorting genetic algorithm II (NSGA-II) \citep{deb2002}, a widely used evolutionary algorithm, was employed to identify optimal designs. The optimization problem was formulated to minimize the drag coefficient and maximize the lift coefficient simultaneously:
\begin{equation}
    \min_{\mathbf{x} \in \mathcal{X}} \quad \mathbf{f}(\mathbf{x}) =
    \begin{bmatrix}
        C_d(\mathbf{x}) \\
        -C_l(\mathbf{x})
    \end{bmatrix}
\end{equation}
where $\mathbf{x}$ denotes the vector of design parameters within the feasible design space $\mathcal{X}$. The optimization process employed a population size of 100 and $10^4$ function evaluations. Polynomial mutation and simulated binary crossover were applied with probabilities set to 0.15 and 0.9, respectively.

Given the conflicting nature of the objectives, a Pareto front was constructed. Each solution on the Pareto front represents a non-dominated design where no other feasible configuration simultaneously improves both objectives. The final optimal design was selected from the Pareto set using the utopia point method \citep{miettinen1999}, which identifies the solution closest to the theoretical ideal (utopia) point defined by maximum $C_l$ and minimum $C_d$ values.

\subsection{Adjoint-based design}

\begin{figure}
    \centering
    \includegraphics[width=\linewidth]{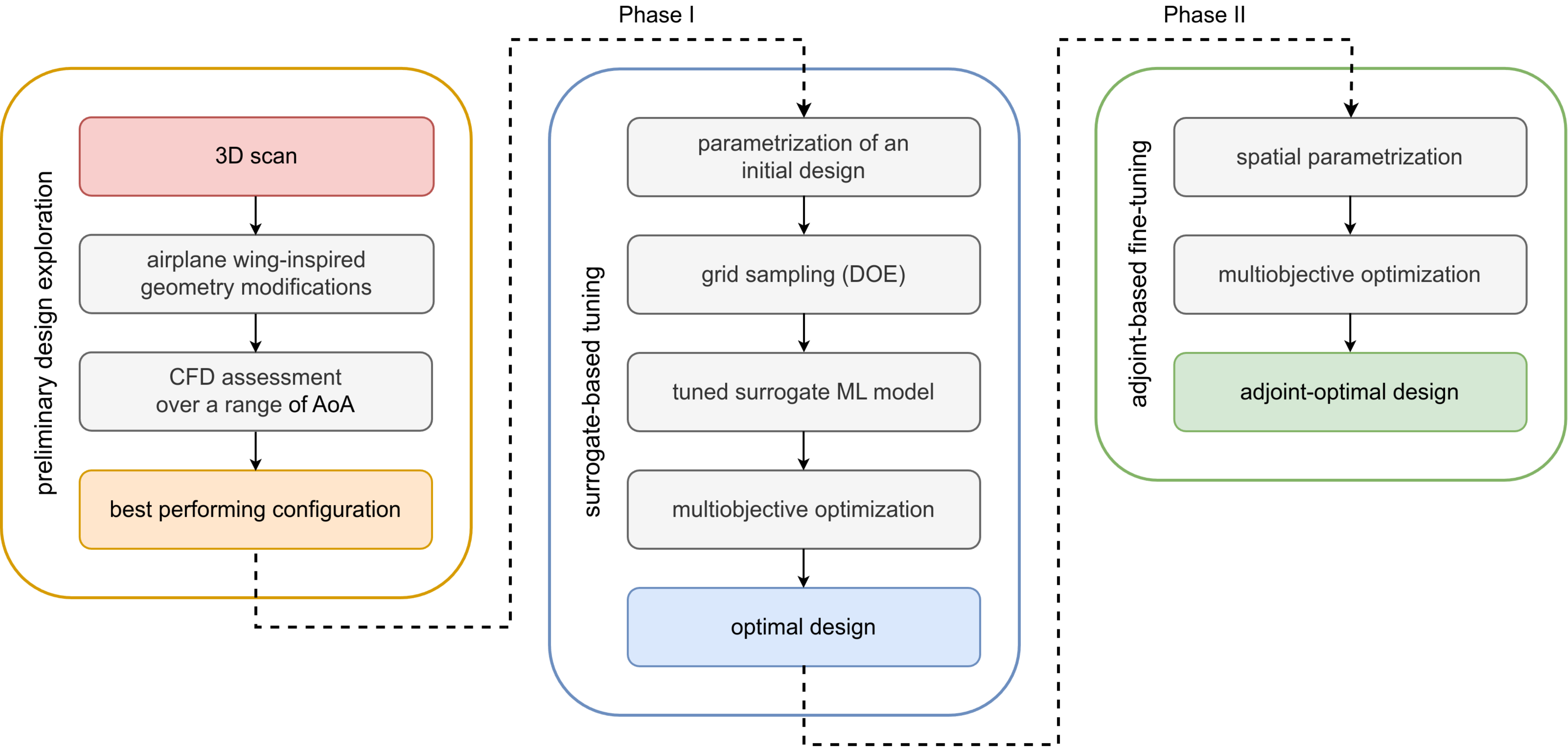}
    \caption{Schematic overview of the ML surrogate-based and adjoint-based optimization methodology.}
    \label{fig:flowchart}
\end{figure}

Following the initial design exploration using a machine-learning-based (ML) surrogate and multi-objective optimization, a gradient-based adjoint optimization approach \citep{giles2000} was employed to further refine the hydrofoil design. An overview of the methodology is shown in Figure \ref{fig:flowchart}. The adjoint method provides exact sensitivity information for all design parameters simultaneously, with a computational cost that is independent of the number of design variables, which makes it suitable for high-dimensional optimization problems.

The adjoint optimization process relies on free-form deformation, i.e., the shape of the hydrofoil is deformed by manipulating the positions of control points. B-spline basis functions define the relationship between the new positions of the control points and the deformation of every point on the hydrofoil's surface and the surrounding volume mesh. The control points effectively serve as the design variables for the optimization problem.

The adjoint optimization was conducted using the adjointOptimizationFoam \citep{papoutsis2021} tool within the OpenFOAM v2312 framework. This tool computes the gradient of the objective function with respect to design parameters by solving the adjoint equations alongside the primal flow equations. The free-form deformation parameterization is implemented natively within adjointOptimizationFoam through its volumetric B-spline morpher. This built-in mesh motion solver creates a three-dimensional control grid around the geometry, where the coordinates of control points serve as design variables. The morpher uses uniform B-spline or Bezier basis functions to smoothly propagate control point movements to the computational mesh. The corresponding adjoint equations are derived by introducing Lagrange multipliers (adjoint variables) and taking the functional derivative of the augmented objective function. The continuous adjoint equations for steady, incompressible, turbulent flow are:
\begin{equation}
    \nabla \cdot \mathbf{u}_a = \frac{\partial J}{\partial p}
\end{equation}
\begin{equation}
    -(\mathbf{u} \cdot \nabla) \mathbf{u}_a + (\nabla \mathbf{u})^T \cdot \mathbf{u}_a = -\frac{1}{\rho}\nabla p_a + \nu \nabla^2 \mathbf{u}_a + \frac{\partial \boldsymbol{\tau}}{\partial \mathbf{u}} \cdot \mathbf{u}_a + \frac{\partial J}{\partial \mathbf{u}}
\end{equation}
where $\mathbf{u}_a$ and $p_a$ are the adjoint velocity and pressure fields, $\mathbf{u}$ and $p$ are the primal flow velocity and pressure, $\boldsymbol{\tau}$ is the turbulent stress tensor, and the right-hand side source terms are the derivatives of the objective functional $J$ with respect to the primal variables. The adjoint optimization problem was formulated as a weighted sum:
\begin{equation}
    \min_{\mathbf{x} \in \mathcal{X}} \quad J = \mathbf{f}(\mathbf{x}) = 1.0 \cdot C_d(\mathbf{x}) - 0.5 \cdot C_l(\mathbf{x})
\end{equation}
subject to the volume constraint:
\begin{equation}
    0.999 \leq \frac{V(\mathbf{x})}{V_0} \leq 1.001
\end{equation}
where $\mathbf{x}$ represents the design variables, $\mathcal{X}$ is the feasible design space, $V(\mathbf{x})$ is the hydrofoil volume, and $V_0$ is the initial volume. The adjoint optimization was run for 500 iterations. The weights were set to 1.0 for the drag coefficient and 0.5 for the lift coefficient based on the analysis of the Pareto front from preliminary multi-objective optimization studies. This analysis revealed that lift increased at a faster rate than drag decreased, indicating a trade-off imbalance. Assigning a lower weight to lift reduced its relative influence in the objective function, helping to balance the competing objectives and prevent the optimizer from disproportionately prioritizing lift gains over drag reduction, which is physically more difficult to achieve. For the conventional multi-objective optimization, this adjustment was not necessary since the complete Pareto front was explicitly computed.

\section{Results}

\subsection{Lift and drag for preliminary configurations}

\begin{figure}[H]
    \centering
    \begin{subfigure}[b]{0.46\linewidth}
        \centering
        \includegraphics[width=\linewidth]{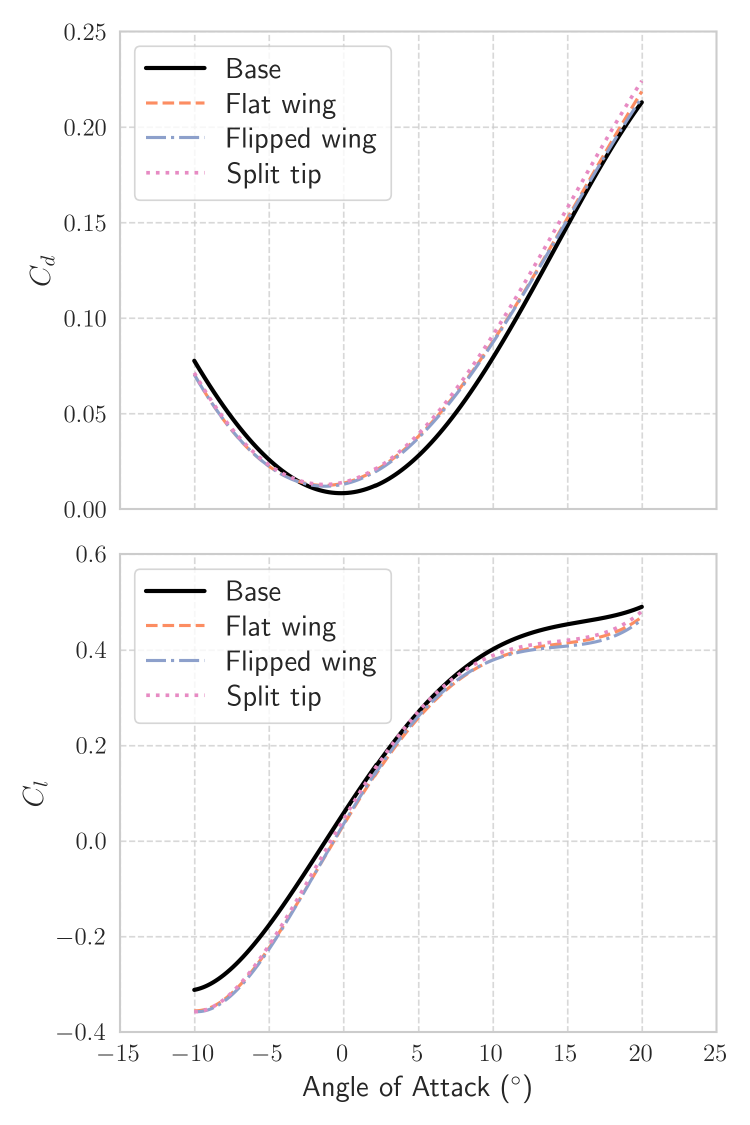}
    \end{subfigure}
    \hfill
    \begin{subfigure}[b]{0.46\linewidth}
        \hspace{-0.5cm}
        \includegraphics[width=\linewidth]{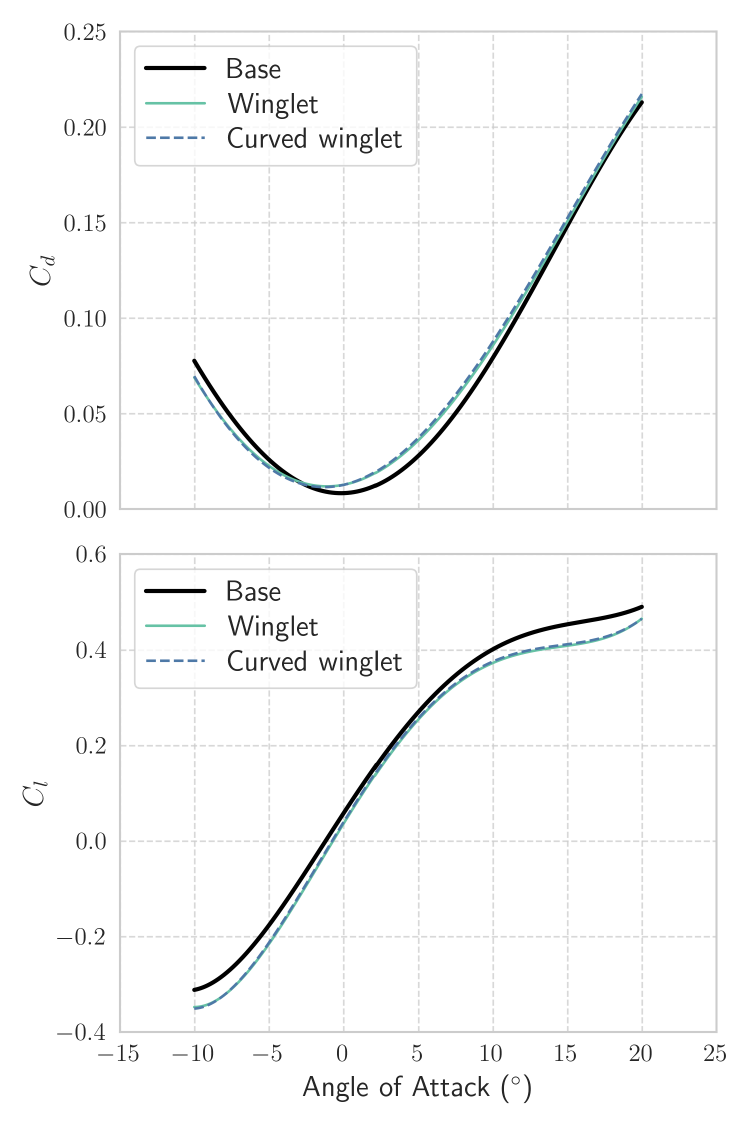}
    \end{subfigure}
    \caption{Drag coefficient $C_d$ (top) and lift coefficient $C_l$ (bottom) across a range of angles of attack for the preliminary front wing designs.}
    \label{fig:designs_results_Cl_Cd}
\end{figure}

The preliminary configurations were assessed by comparing the drag coefficient $C_d$ and lift coefficient $C_l$ against the baseline design. The results, presented as a function of the angle of attack, are shown in Figure \ref{fig:designs_results_Cl_Cd}.

The baseline design outperforms all other configurations across the tested AoA range, which is expected given that the preliminary designs were chosen arbitrarily. For AoA between approximately -10$^{\circ}$ and -2$^{\circ}$, all novel configurations achieve slightly lower $C_d$ than the baseline, though the differences are minimal. Between -2$^{\circ}$ and 15$^{\circ}$ AoA, the baseline design is clearly superior. Beyond 15$^{\circ}$ AoA, the performance gap narrows. Among the novel configurations, winglet variants exhibit the lowest drag.

All novel configurations consistently generate lower lift than the baseline hydrofoil across the tested AoA range. The difference is minor between 0$^{\circ}$ and 5$^{\circ}$ AoA, but becomes more pronounced outside this range. Within this range, the split-tip design shows the highest lift among the preliminary configurations.

\begin{figure}[H]
    \centering
    \includegraphics[width=0.9\linewidth]{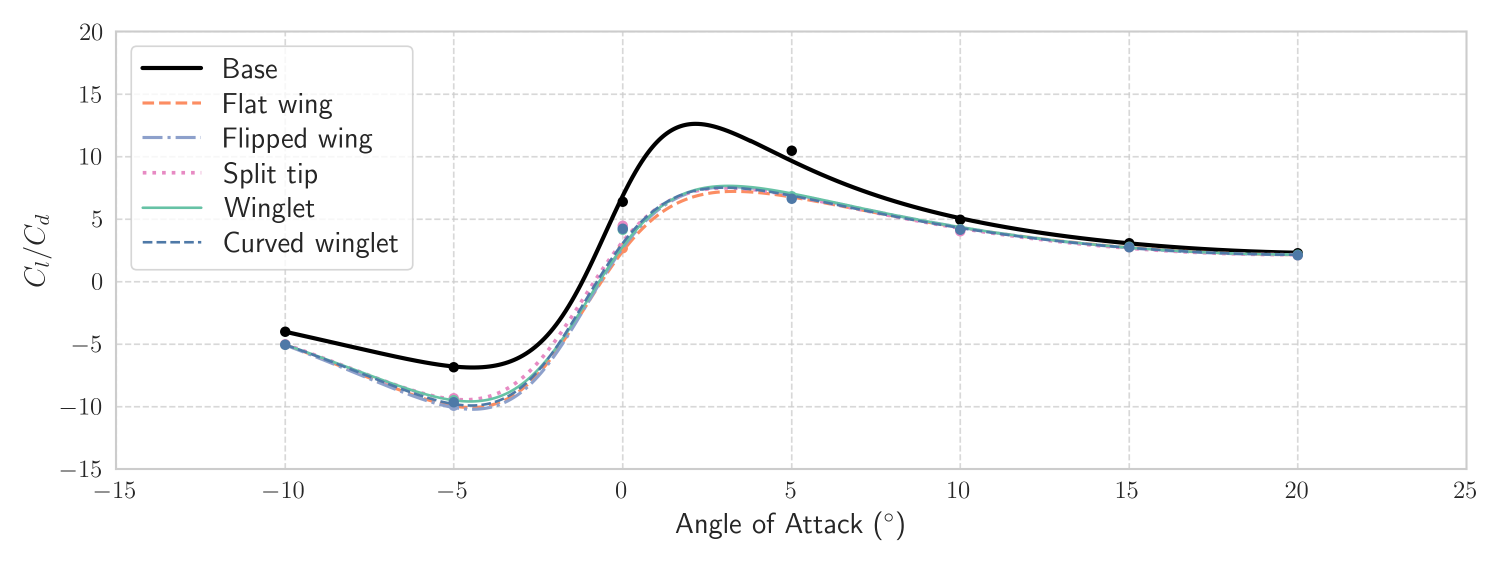}
    \caption{Hydrodynamic efficiency $C_l / C_d$ across a range of angles of attack for the preliminary front wing designs.}
    \label{fig:designs_results_Cl_Cd_ratio}
\end{figure}

The lift-to-drag ratio $C_l/C_d$, shown in Figure \ref{fig:designs_results_Cl_Cd_ratio}, is used as the primary metric to select the most promising preliminary configuration. This ratio quantifies the efficiency of the hydrofoil, i.e., the amount of lift generated for a given amount of drag. A higher ratio indicates better efficiency. Similar to the $C_d$ and $C_l$ trends, the baseline design achieves the highest efficiency. The performance gap is most pronounced at 5$^{\circ}$ AoA, where the baseline design's efficiency exceeds the best preliminary configurations by approximately 40\%. Between -2$^{\circ}$ and 2$^{\circ}$ AoA, the difference is smaller, up to approximately 15\%.

Based on overall performance, the winglet configuration was selected for further assessment. Although its peak performance is lower than the baseline, it outperforms the remaining preliminary configurations across the AoA range. The general performance trend can be attributed to the fact that the geometries were created in a rough, inspiration-driven manner from existing wing designs, without adhering to specific design guidelines or an optimization strategy.

\subsection{Surrogate-based multi-objective optimization}

\begin{figure}[H]
	\centering
	\includegraphics[width=0.9\linewidth]{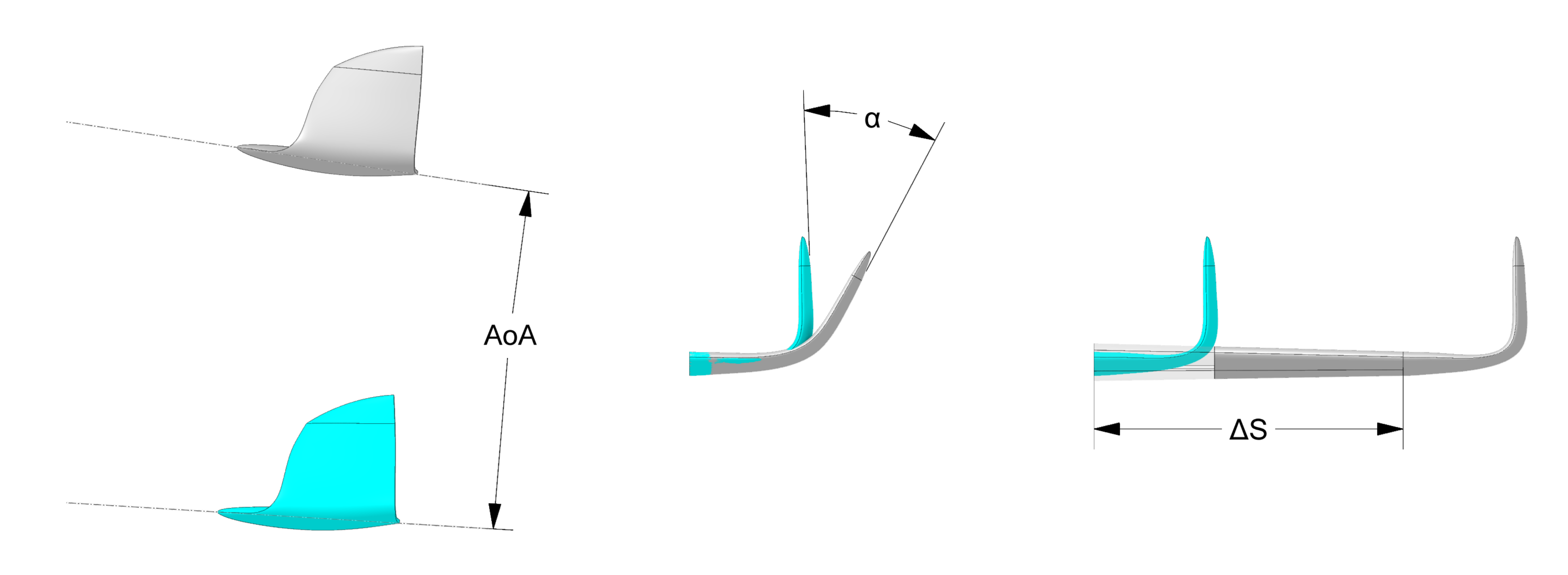}
	\caption{DOE design parameters: angle of attack, tip angle, and span.}
	\label{fig:doe}
\end{figure}

The spanwise variation in the front wing's cross-sections has considerable geometric complexity, which consequently limited the parametric study to two principal design variables. A comprehensive parametrization of the geometry would entail a high-dimensional variable space, which was deemed impractical. The selected parameters were the winglet tip angle $\alpha$ and wing span length $\Delta$S, with the angle of attack as an operational parameter. The chosen parameters are shown in Figure \ref{fig:doe}.

The winglet tip angle was defined relative to the vertical axis of the winglet, with the baseline design corresponding to $\alpha = 0^{\circ}$ (a vertically oriented winglet). Span length modifications were implemented by offsetting the winglet and lofting between the root and tip cross-sections. The investigated span changes included reductions and extensions of the front wing. The AoA range spans a wide range of operational conditions.

A full-factorial DOE approach was employed, investigating all possible combinations of the three parameters. This resulted in a total of 72 CFD simulations. The complete parameter matrix is summarized in Table \ref{tab:DOE_parameters}.

\begin{table}[H]
    \centering
    \footnotesize
    \caption{Design parameters and values considered in the DOE approach.}
    \label{tab:DOE_parameters}
    \begin{tabular}{@{}lccc@{}}
        \hline
        \textbf{Design parameter} & \textbf{Symbol} & \textbf{Unit} & \textbf{Values} \\
        \hline
        Angle of attack & AoA & $^{\circ}$ & -5,\ -2.5,\ 0,\ 2.5,\ 5,\ 10 \\
        Tip angle & $\alpha$ & $^{\circ}$ & 0,\ 30,\ 60 \\
        Span & $\Delta$S & mm & -30,\ 0,\ 30,\ 80 \\
        \hline
    \end{tabular}
\end{table}

\begin{figure}[H]
	\centering
	\includegraphics[width=0.6\linewidth]{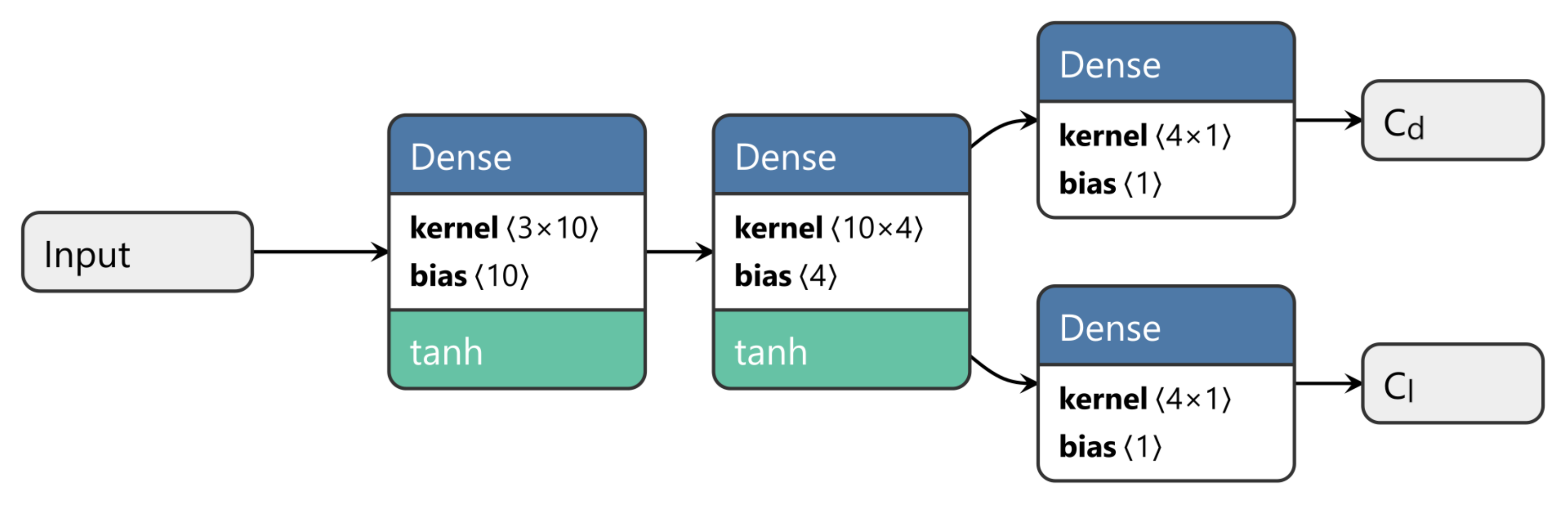}
	\caption{Schematic of the MLP architecture.}
	\label{fig:mlp_model}
\end{figure}

Data generated in the DOE approach was used to train an ML surrogate according to the methodology outlined in Section \ref{sec:data_driven_design}. The tuned MLP uses three inputs (design parameters) and two fully connected hidden layers with 10 and 4 neurons, respectively. The network architecture branches into two separate output layers, each producing a single scalar value corresponding to drag and lift coefficients, respectively. The tuned surrogate achieved a combined $R^2$ of 0.97, with an RRMSE of 4.52\%. A schematic representation of the MLP architecture is given in Figure \ref{fig:mlp_model}.

\begin{figure}[H]
    \centering
    \begin{subfigure}[b]{0.45\linewidth}
        \centering
        \includegraphics[width=\linewidth]{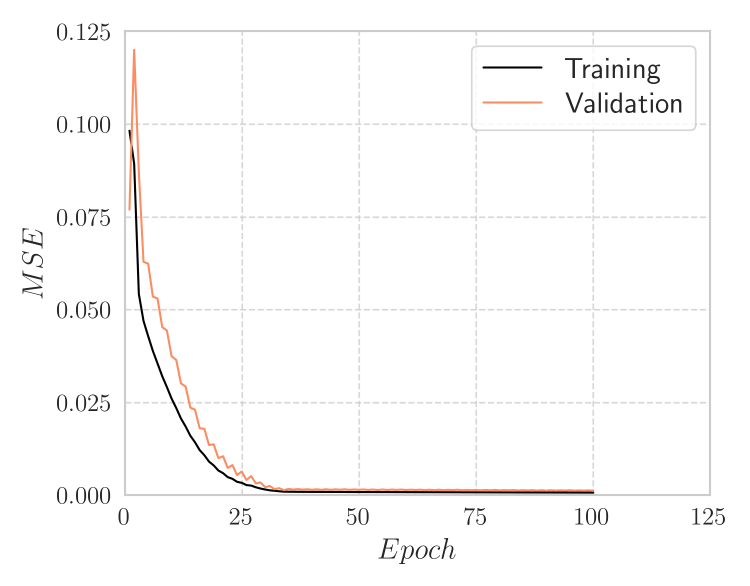}
        \caption{}
    \end{subfigure}
    \hfill
    \begin{subfigure}[b]{0.45\linewidth}
        \centering
        \includegraphics[width=\linewidth]{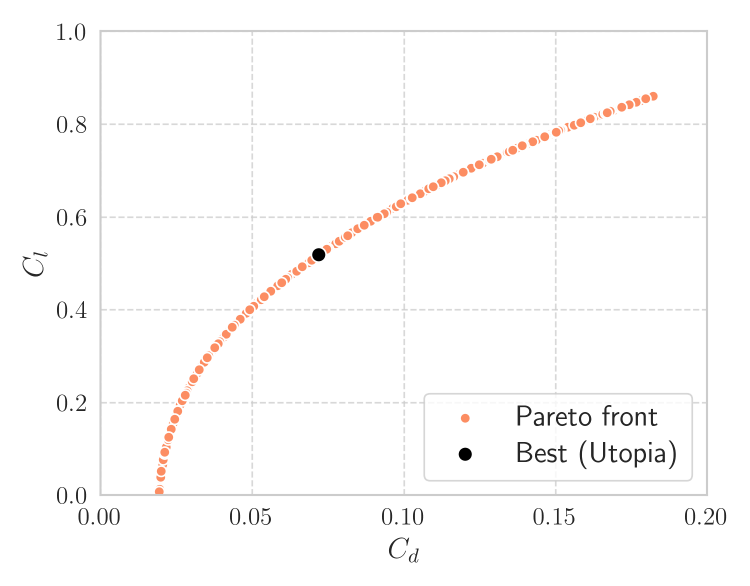}
        \caption{}
    \end{subfigure}
    \caption{Surrogate training and optimization results: (a) training and validation loss curves of the tuned MLP surrogate model, (b) Pareto front resulting from the multi-objective optimization, where the trained surrogate model is used as the objective function.}
    \label{fig:loss_curve_pareto_front}
\end{figure}

The tuned MLP surrogate was coupled with the NSGA-II algorithm to identify the optimal hydrofoil design. The optimum was selected based on the utopia point approach. Figure \ref{fig:loss_curve_pareto_front} shows the training and validation loss curves of the MLP as well as the obtained Pareto front and optimal result. The Pareto front indicates that the drag coefficient, $C_d$, is typically lower, approximately one-fifth of the lift coefficient, $C_l$. This implies that in a weighted multi-objective approach, equally weighting both objectives will naturally emphasize lift improvements over drag reduction.

\begin{figure}[H]
    \centering
    \begin{subfigure}[b]{0.495\linewidth}
        \centering
        \includegraphics[trim={10cm 10cm 7.5cm 10cm}, clip, width=\linewidth]{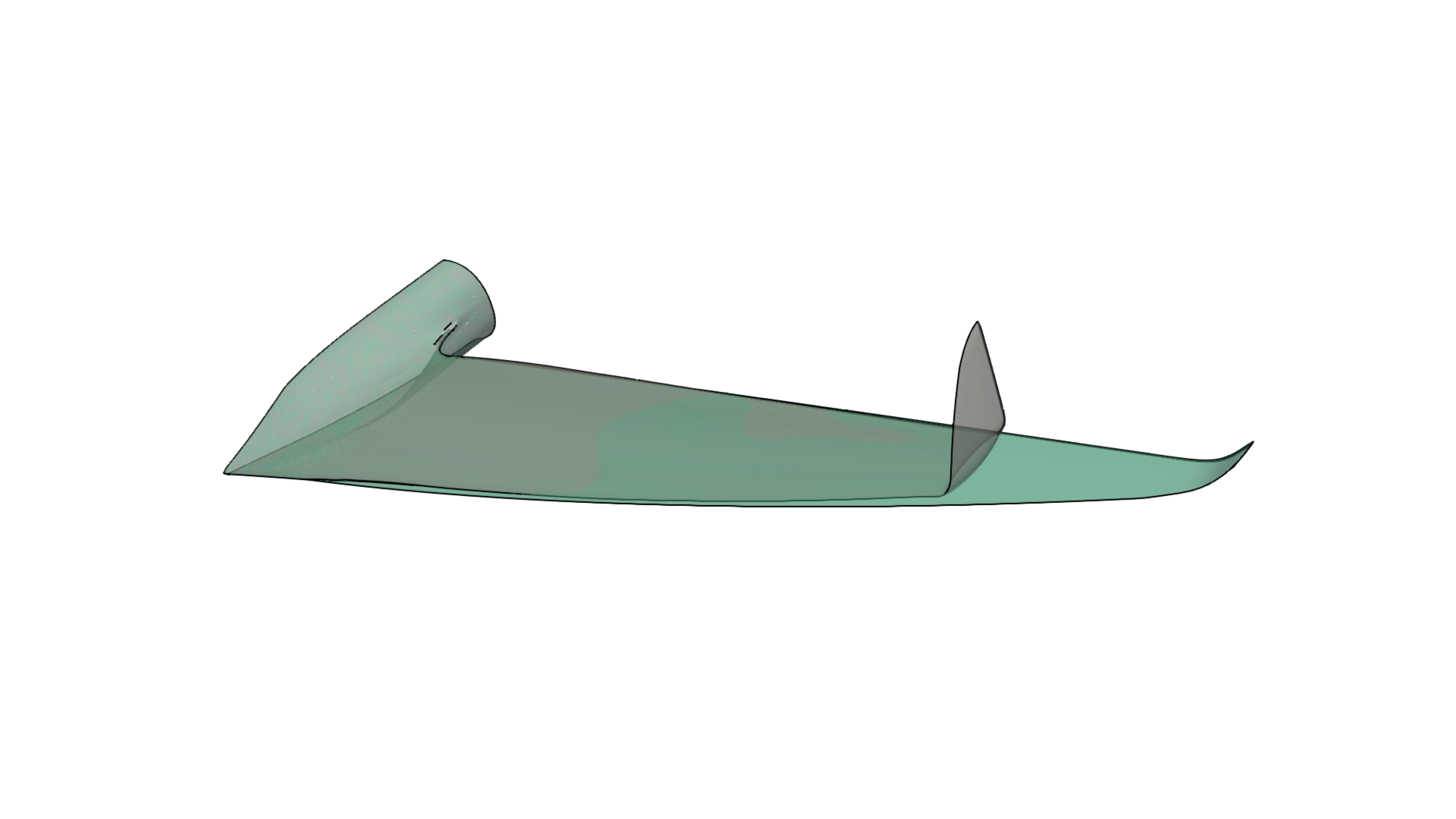}
    \end{subfigure}
    \caption{Comparison of the initial winglet geometry (gray) and the best design obtained from the ML-MO optimization workflow (green).}
    \label{fig:doe_results}
\end{figure}

Figure \ref{fig:doe_results} shows the baseline winglet geometry and the optimized shape obtained through the ML-MO optimization approach for 0$^{\circ}$ AoA. At 0$^{\circ}$ AoA, the ML-MO optimization yields an optimum with a winglet tip angle of $\alpha = 52^\circ$ and a span length extension of $\Delta$S$ = 28$ mm. The optimal design maintains the essential winglet configuration, but with notable changes: the vertical component of the winglet is reduced and bent downward, effectively extending the wing's span. In planform view, the forward wing section is enlarged, while maintaining approximately the same frontal projected area as the initial design.

\begin{figure}[H]
	\centering
	\includegraphics[width=\linewidth]{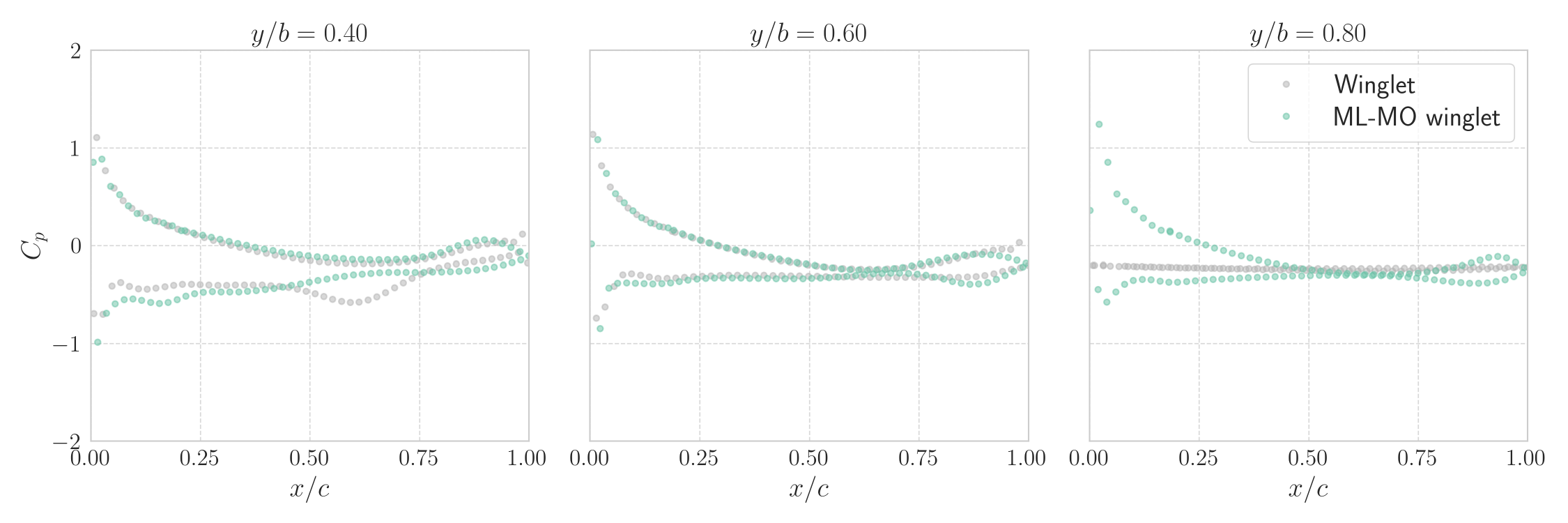}
	\caption{Pressure coefficient $C_p$ at spanwise sections $y/b=0.4$, $0.6$, and $0.8$ (left to right) for the initial winglet and ML-MO design.}
	\label{fig:winglet_results_Cp}
\end{figure}

Pressure coefficient distributions at different spanwise sections are shown in Figure \ref{fig:winglet_results_Cp}. Compared to the initial winglet configuration, the ML-MO design consistently exhibits a greater pressure differential between the suction and pressure sides, indicating higher sectional lift. A minor irregularity is visible at $y/b = 0.4$ near $x/c = 0.6$. Further along the span, the observed trends reflect the geometric differences between the two designs. The ML-MO design is more laid out rather than angled upward like the initial configuration. This change in orientation explains the altered pressure distribution toward the tip region.

\begin{figure}[H]
    \centering
    \includegraphics[width=0.9\linewidth]{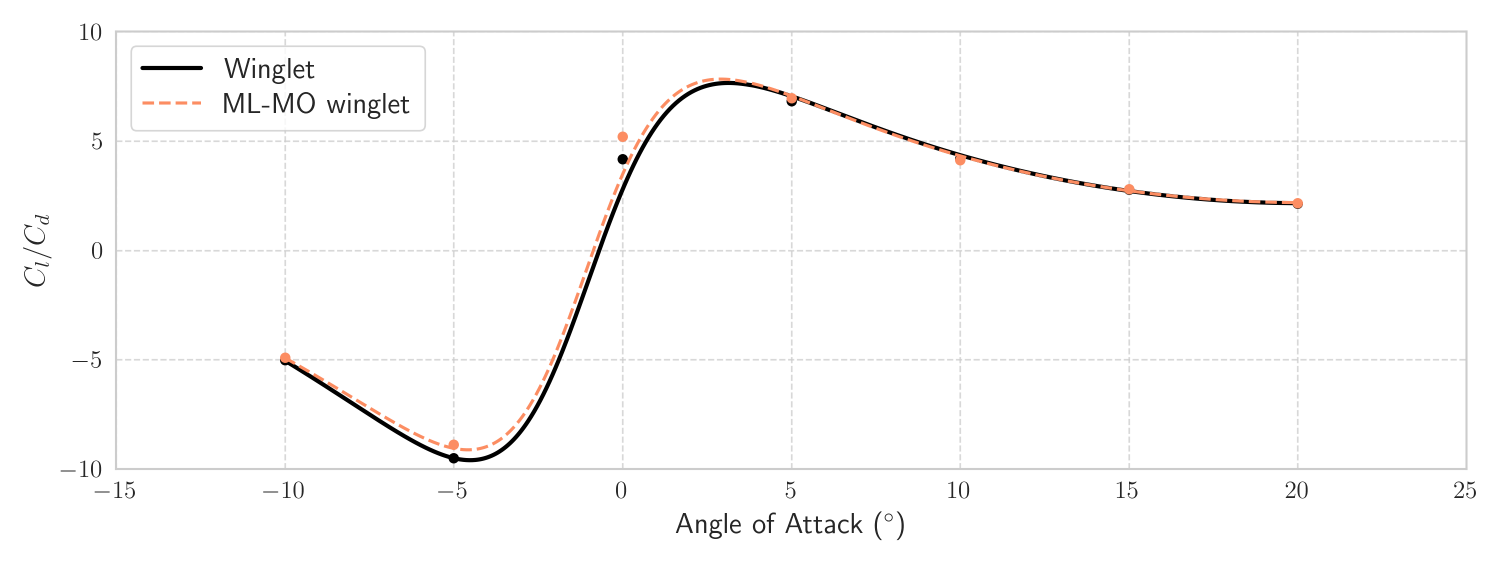}
    \caption{Hydrodynamic efficiency $C_l / C_d$ across a range of angles of attack for the base winglet and ML-MO-optimized configuration.}
    \label{fig:winglet_results_Cl_Cd_ratio}
\end{figure}

Based on the results shown in Figure \ref{fig:winglet_results_Cl_Cd_ratio}, the optimized design provides moderate performance improvements. At angles of attack between -5$^{\circ}$ and 5$^{\circ}$, efficiency improves, with the largest gain of about 24\% at 0$^{\circ}$ AoA. The limited performance improvement can be attributed to two main factors. First, the multi-objective optimization was conducted at a specific AoA, resulting in the most significant improvement being localized to that operating condition. Second, due to the limited number of design parameters, the scope of design changes was constrained, limiting performance improvements.

\subsection{Adjoint optimization}

\begin{figure}[H]
    \centering
    \begin{subfigure}[b]{\linewidth}
        \centering
        \includegraphics[trim={4cm 8cm 2cm 6cm}, clip, width=0.75\linewidth]{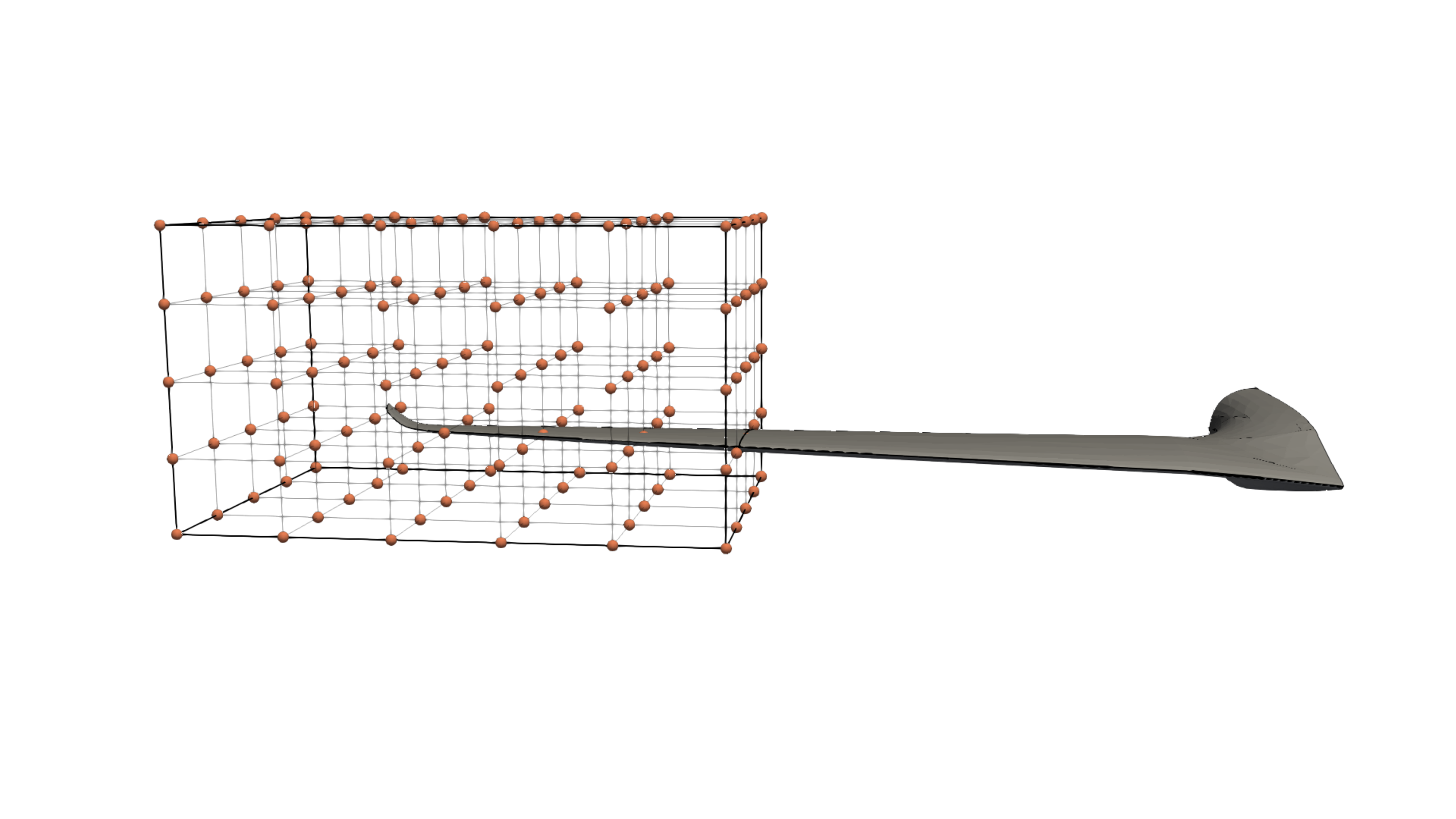}
        \caption{}
    \end{subfigure}
    \vfill
    \begin{subfigure}[b]{\linewidth}
        \centering
        \begin{subfigure}[b]{0.495\linewidth}
            \centering
            \includegraphics[width=\linewidth]{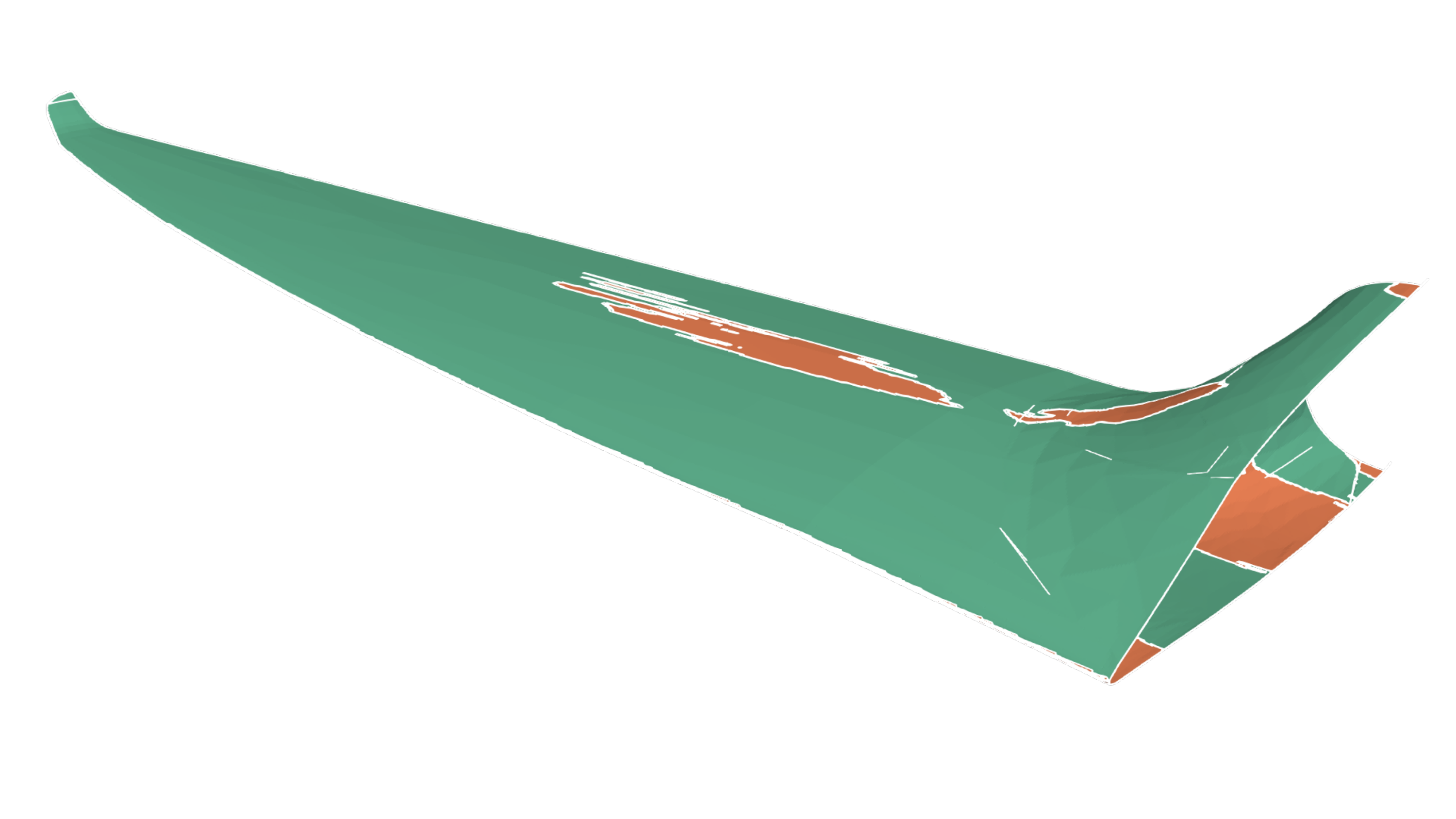}
            \caption{}
        \end{subfigure}
        \hfill
        \begin{subfigure}[b]{0.495\linewidth}
            \centering
            \includegraphics[width=\linewidth]{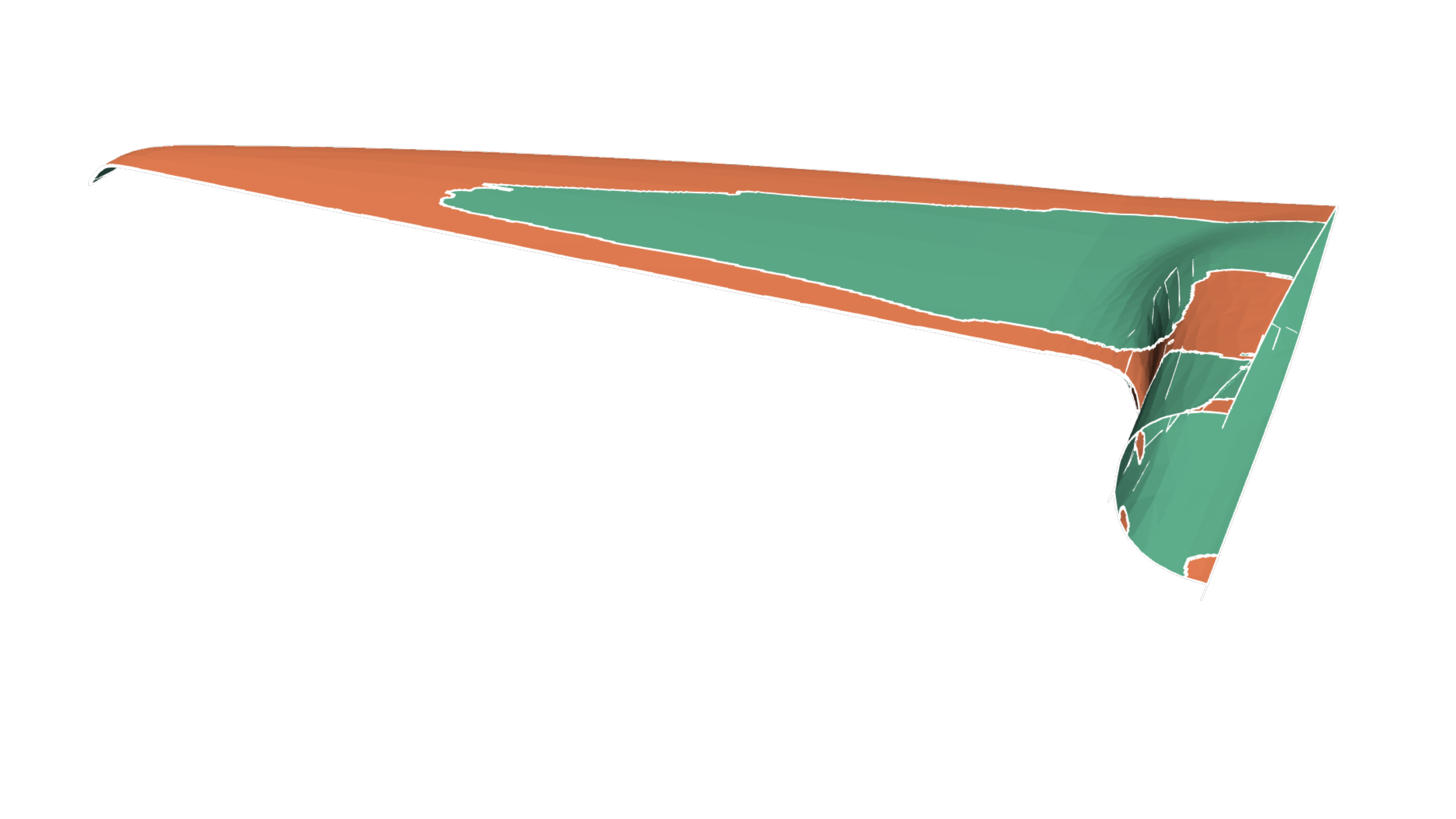}
            \caption{}
        \end{subfigure}
    \end{subfigure}
    \caption{Adjoint optimization of the ML-MO improved winglet: (a) control points defining the deformable domain, (b) top view of smoothed adjoint sensitivities, with green indicating inward and orange indicating outward deformation, (c) bottom view of smoothed adjoint sensitivities, employing the same color convention.}
    \label{fig:doe_winglet_adjoint}
\end{figure}

Given the modest performance improvements obtained through the ML-MO optimization approach, a logical next step is to allow greater flexibility in the design. Adjoint optimization is a natural option, as it allows significant flexibility and efficient simultaneous tuning. The optimal ML-MO design at 0$^{\circ}$ AoA was selected as the basis for adjoint optimization. A deformable domain was defined using 150 control points (arranged as 6 points along the span and 5 points in each of the other two directions). To ensure robust optimization, adjointOptimisationFoam was extended with an automatic remeshing procedure that is triggered whenever the mesh fails the in-built quality check. Afterward, the process continues with the next iteration. The positioning of the control points and the extents of the deformable domain were determined based on adjoint sensitivities. Figure \ref{fig:doe_winglet_adjoint} shows the control points and adjoint sensitivities for the winglet design.

\begin{figure}[H]
    \centering
    \includegraphics[width=0.9\linewidth]{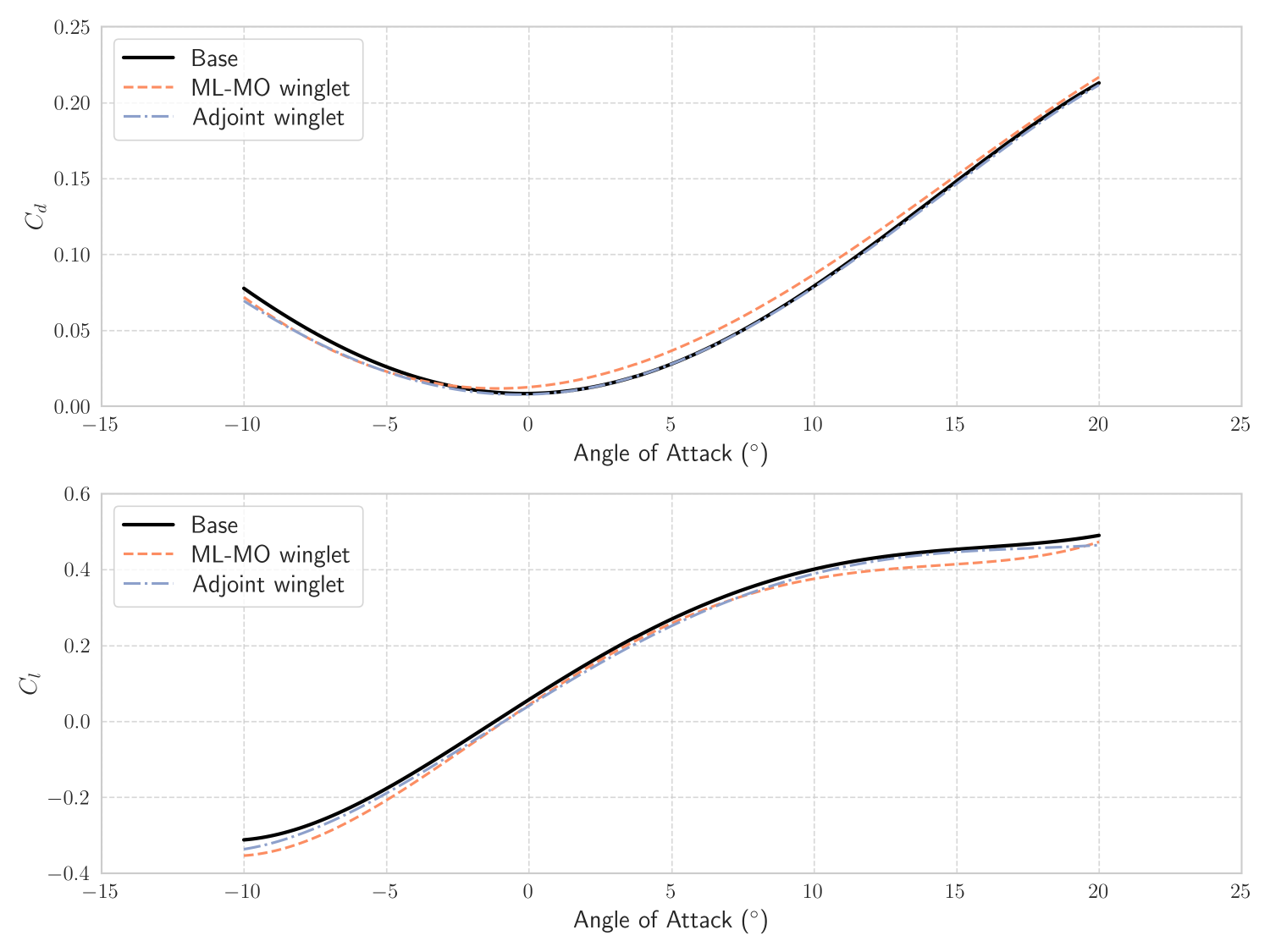}
    \caption{Drag coefficient $C_d$ (top) and lift coefficient $C_l$ (bottom) across a range of angles of attack for the base front wing design, the ML-MO-optimized winglet design, and the adjoint-improved winglet design.}
    \label{fig:adjoint_doe_winglet_results_Cd_Cl}
\end{figure}

Figure \ref{fig:adjoint_doe_winglet_results_Cd_Cl} shows the drag coefficient and lift coefficient for the baseline, ML-MO-optimized, and adjoint-optimized designs. As can be observed, $C_d$ for the adjoint design is consistently on par and, below 0$^{\circ}$ AoA, outperforms the baseline design. Compared to the ML-MO design, the drag coefficient is lower across the entire AoA range. In terms of $C_l$, the observed trend differs. The adjoint-improved design clearly outperforms the ML-MO design below 0$^\circ$ and above 10$^\circ$ AoA; however, between 0$^\circ$ and 10$^\circ$ AoA, lift is slightly lower. When compared to the baseline design, similar conclusions can be drawn: trends are well matched (though values are lower), with a slight drop-off towards the upper end of the AoA range.

\begin{figure}[H]
    \centering
    \includegraphics[width=0.9\linewidth]{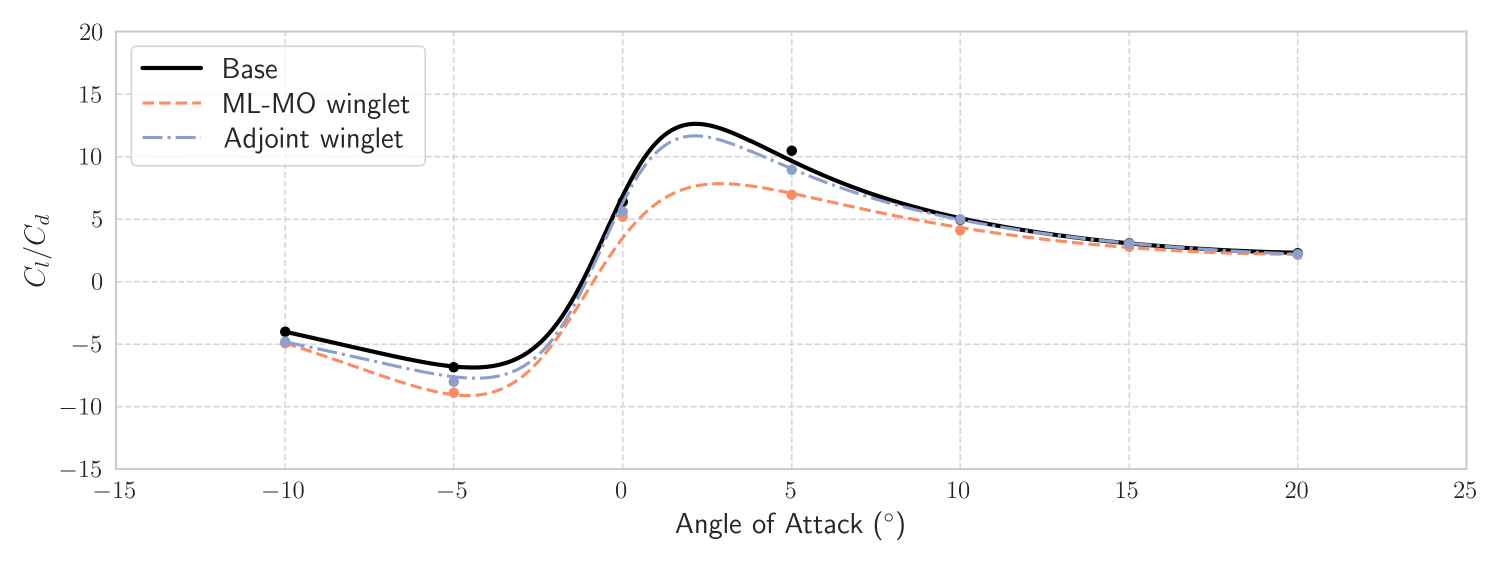}
    \caption{Hydrodynamic efficiency $C_l / C_d$ across a range of angles of attack for the base front wing design, the ML-MO-optimized winglet design, and the adjoint-improved winglet design.}
    \label{fig:adjoint_doe_winglet_results_Cl_Cd_ratio}
\end{figure}

The results for the $C_l / C_d$ ratio, shown in Figure \ref{fig:adjoint_doe_winglet_results_Cl_Cd_ratio}, confirm previous observations. Despite the reduced $C_d$, the adjoint-improved design is unable to match the baseline's performance due to its lower $C_l$. Around 0$^\circ$ AoA, the differences are up to 12\% at worst. Compared to the ML-MO best design, however, the new design represents a clear improvement across the entire AoA range.

\begin{figure}[H]
    \centering
    \begin{subfigure}[b]{\linewidth}
        \centering
        \includegraphics[trim={4cm 4cm 2cm 10cm}, clip, width=0.75\linewidth]{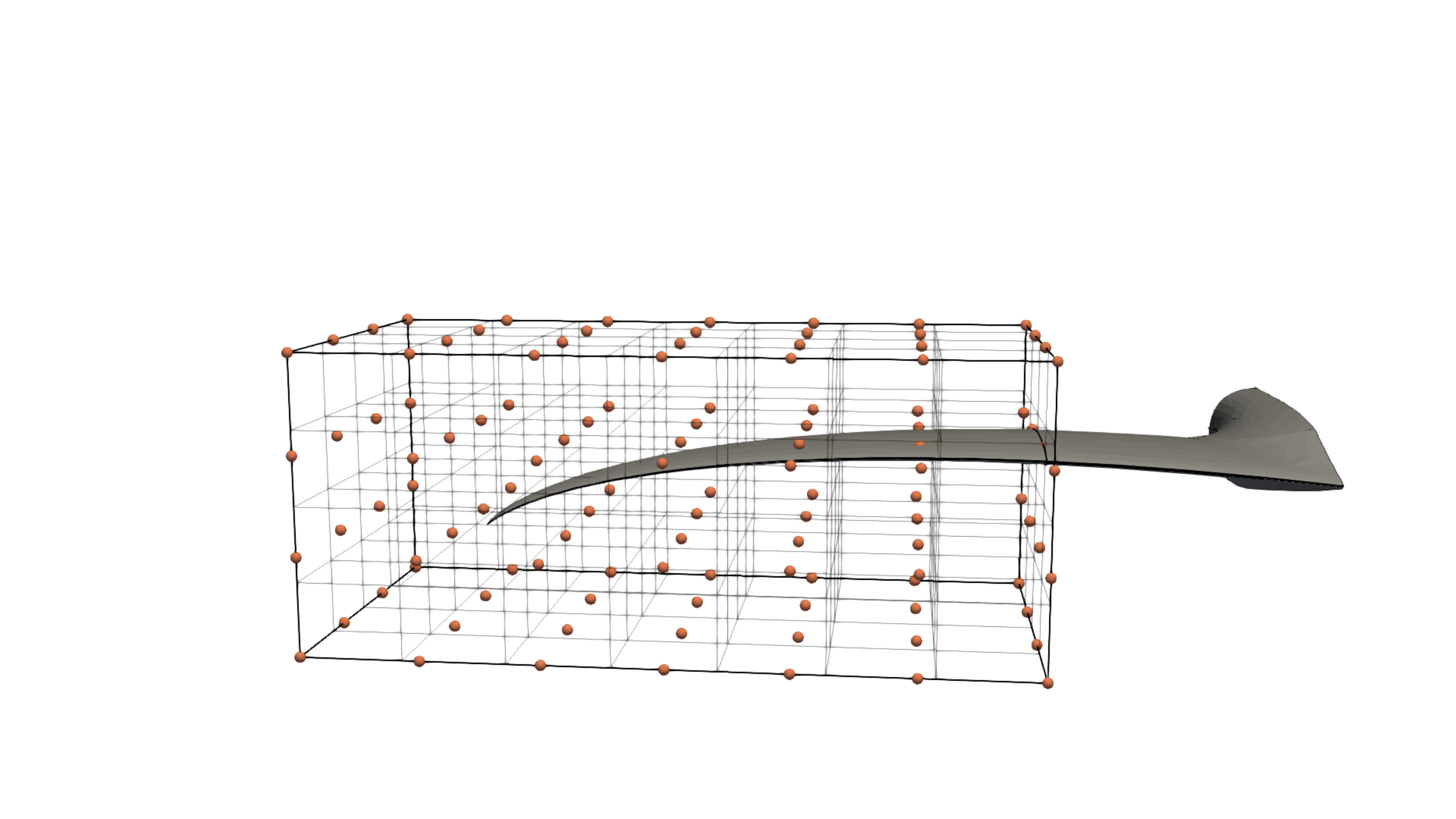}
        \caption{}
    \end{subfigure}
    \vfill
    \begin{subfigure}[b]{\linewidth}
        \centering
        \begin{subfigure}[b]{0.495\linewidth}
            \centering
            \includegraphics[width=\linewidth]{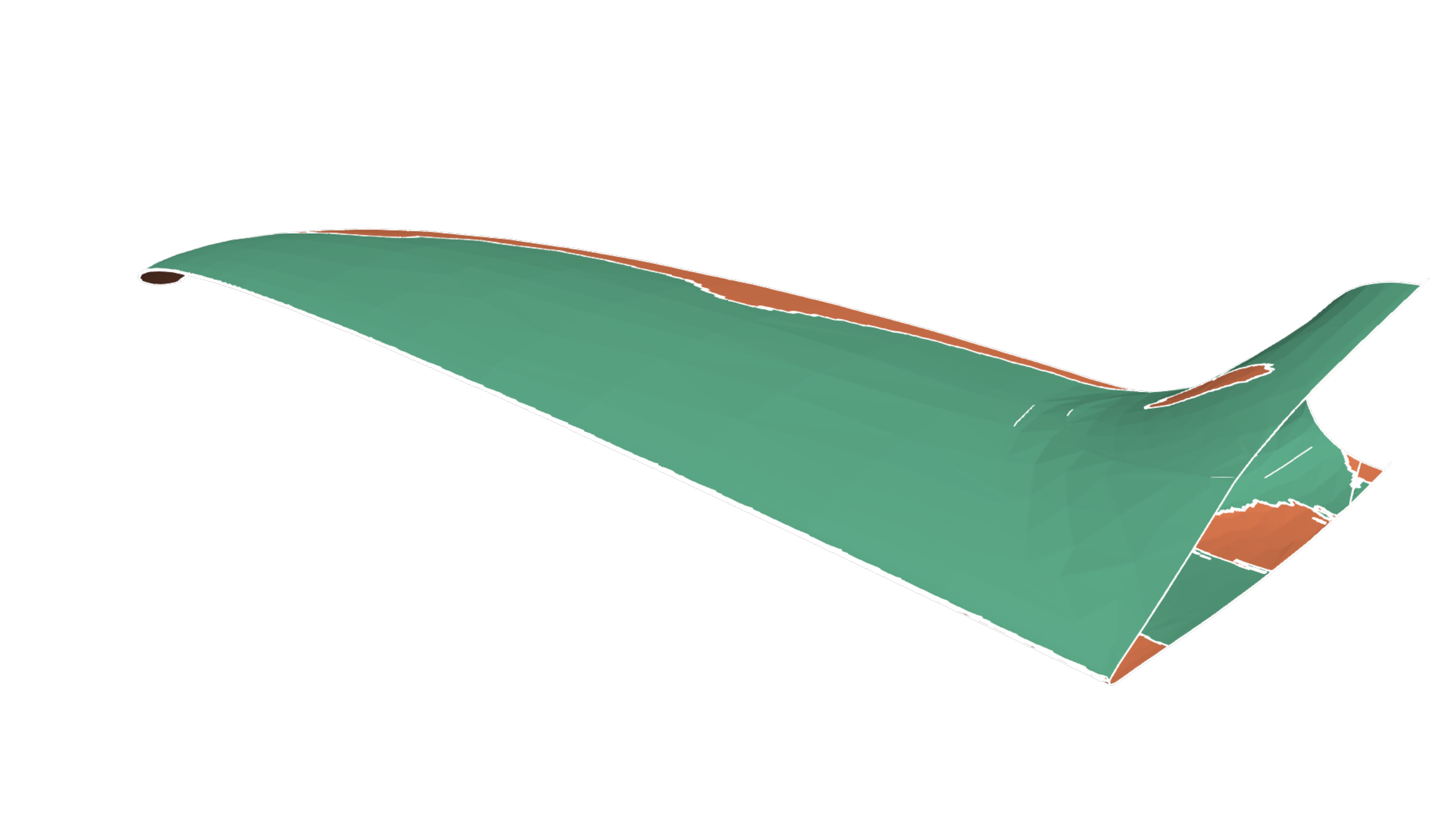}
            \caption{}
        \end{subfigure}
        \hfill
        \begin{subfigure}[b]{0.495\linewidth}
            \centering
            \includegraphics[width=\linewidth]{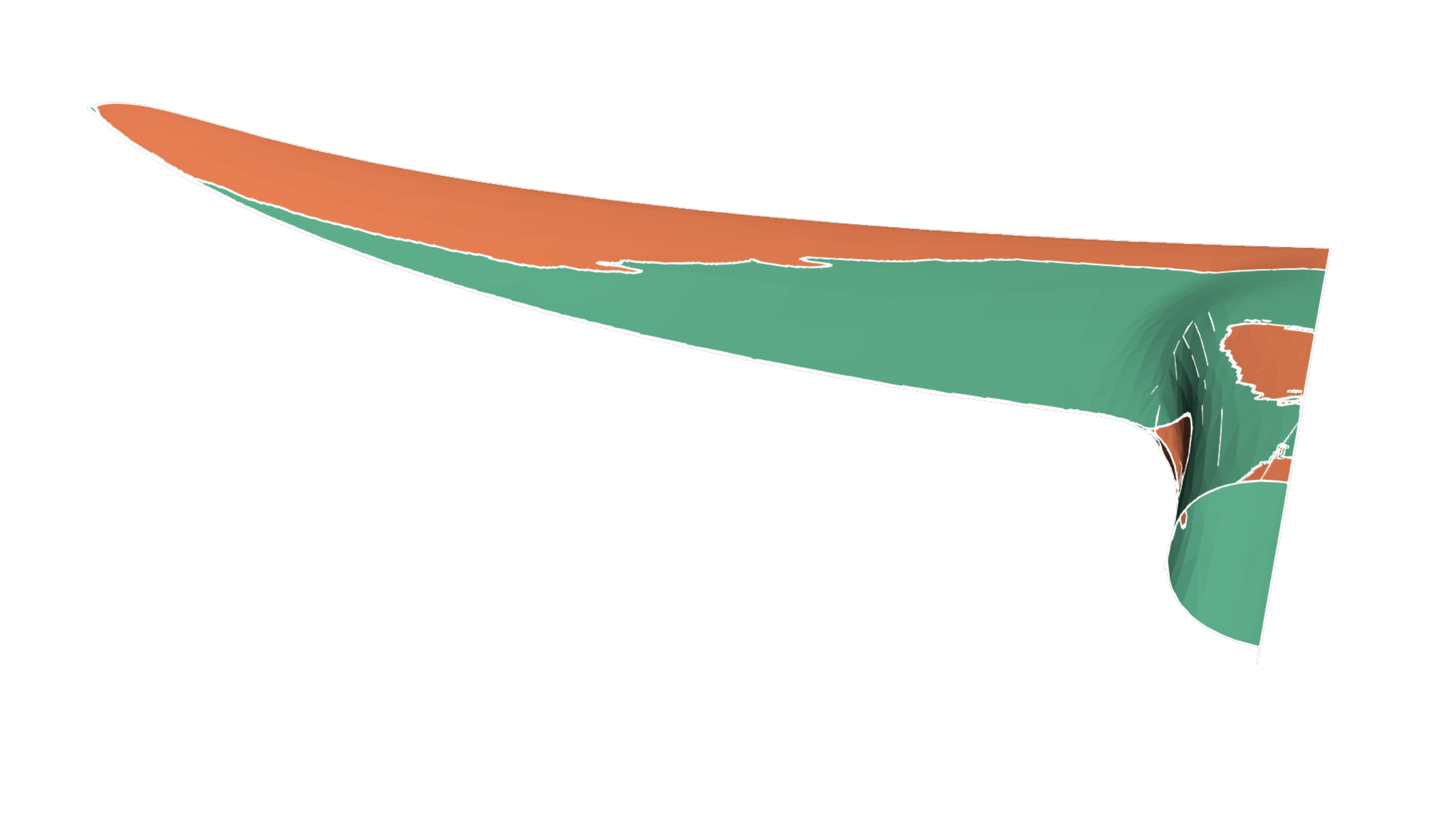}
            \caption{}
        \end{subfigure}
    \end{subfigure}
    \caption{Adjoint optimization of the base wing: (a) control points defining the deformable domain, (b) top view of smoothed adjoint sensitivities, with green indicating inward and orange indicating outward deformation, (c) bottom view of smoothed adjoint sensitivities, employing the same color convention.}
    \label{fig:base_adjoint}
\end{figure}

Despite quantifiable improvements achieved with the adjoint-improved winglet design, it is evident that the baseline design remains superior. The underlying winglet concept appears to be inherently limiting in terms of overall performance improvement. While the winglet can, depending on the operating conditions, achieve performance comparable to the baseline, it is an inferior design solution. To find a better solution, the baseline design was subsequently subjected to adjoint optimization as well. Analogous to the ML-MO-optimized design, the deformable domain was defined based on computed sensitivities. In total, 112 control points were defined (7 points along the span and 4 points in each of the remaining two directions). The setup is shown in Figure \ref{fig:base_adjoint}.

\begin{figure}[H]
    \centering
    \includegraphics[width=0.9\linewidth]{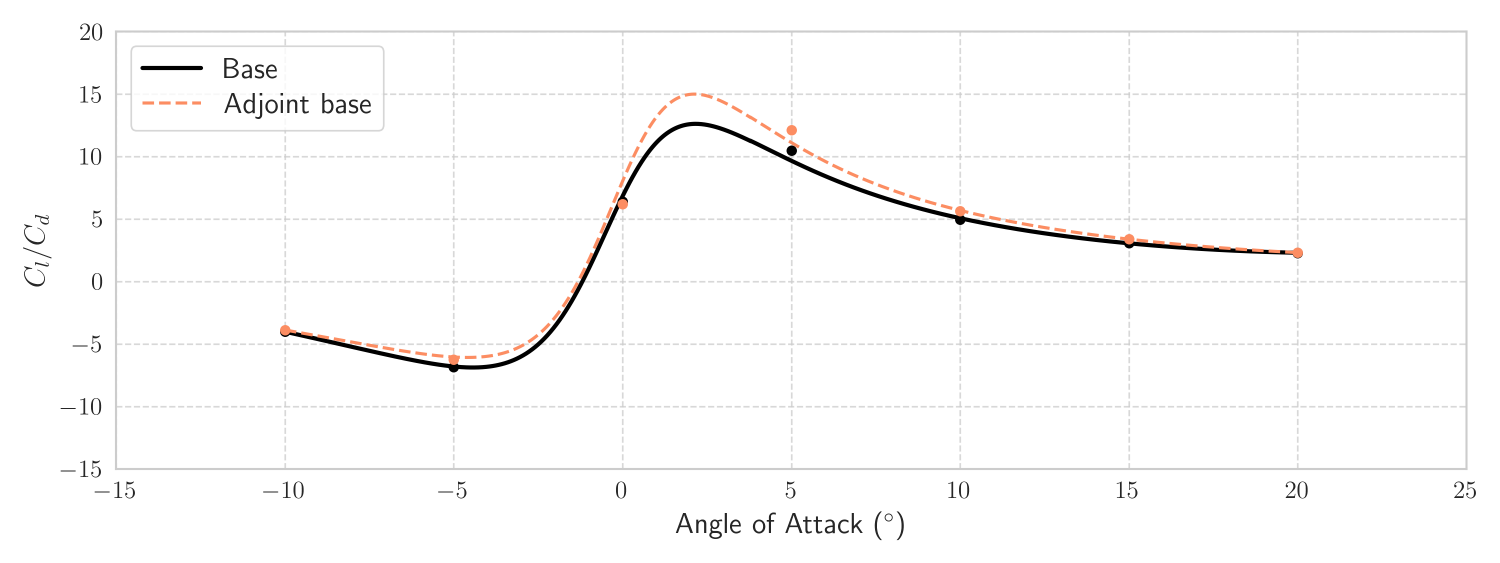}
    \caption{Hydrodynamic efficiency $C_l / C_d$ across a range of angles of attack for the base front wing design and adjoint-improved base design.}
    \label{fig:adjoint_base_results_Cl_Cd_ratio}
\end{figure}

The adjoint-improved baseline design achieves similar or improved $C_l/C_d$ performance compared to the original baseline design across the entire AoA range, as shown in Figure \ref{fig:adjoint_base_results_Cl_Cd_ratio}. The adjoint-improved design performs at least on par with the baseline, with substantial improvements at the mid-range AoA. Around $2^\circ$ AoA, the optimized design provides an efficiency improvement of approximately 20\%. However, these gains diminish toward the lower and upper ends of the AoA range.

Table \ref{tab:cl_cd_ratios} provides a quantitative comparison of the lift-to-drag ratios for baseline and improved configurations considered in this study. The adjoint-improved baseline design demonstrates the highest gains. At extreme angles, the differences between configurations become less pronounced, with all designs converging to similar $C_l / C_d$ values.

\begin{table}[H]
	\centering
	\footnotesize
	\caption{Lift-to-drag ratio ($C_l / C_d$) for baseline and optimized configurations at various AoA.}
	\label{tab:cl_cd_ratios}
	\begin{tabular}{@{}lccccccc@{}}
		\hline
		\multirow{2}{*}{\textbf{\boldmath$C_l / C_d$}} & \multicolumn{7}{c}{\textbf{AoA}} \\
		& \textbf{-10$^\circ$} & \textbf{-5$^\circ$} & \textbf{0$^\circ$} & \textbf{5$^\circ$} & \textbf{10$^\circ$} & \textbf{15$^\circ$} & \textbf{20$^\circ$} \\
		\hline
		Base & -4.02 & -6.85 & 6.39 & 10.47 & 4.95 & 3.07 & 2.30 \\
		Winglet & -5.03 & -9.50 & 4.17 & 6.83 & 4.21 & 2.79 & 2.14 \\
		ML-MO winglet & -4.92 & -8.89 & 5.19 & 6.96 & 4.13 & 2.80 & 2.17 \\
		Adjoint winglet & -4.82 & -8.02 & 5.61 & 8.95 & 4.50 & 3.05 & 2.19 \\
		Adjoint base & -3.87 & -6.25 & 5.93 & 12.35 & 5.64 & 3.38 & 2.31 \\
		\hline
	\end{tabular}
\end{table}

\begin{figure}[H]
    \centering
    \begin{subfigure}[b]{0.495\linewidth}
        \centering
        \includegraphics[trim={10cm 12.5cm 7.5cm 10cm}, clip, width=\linewidth]{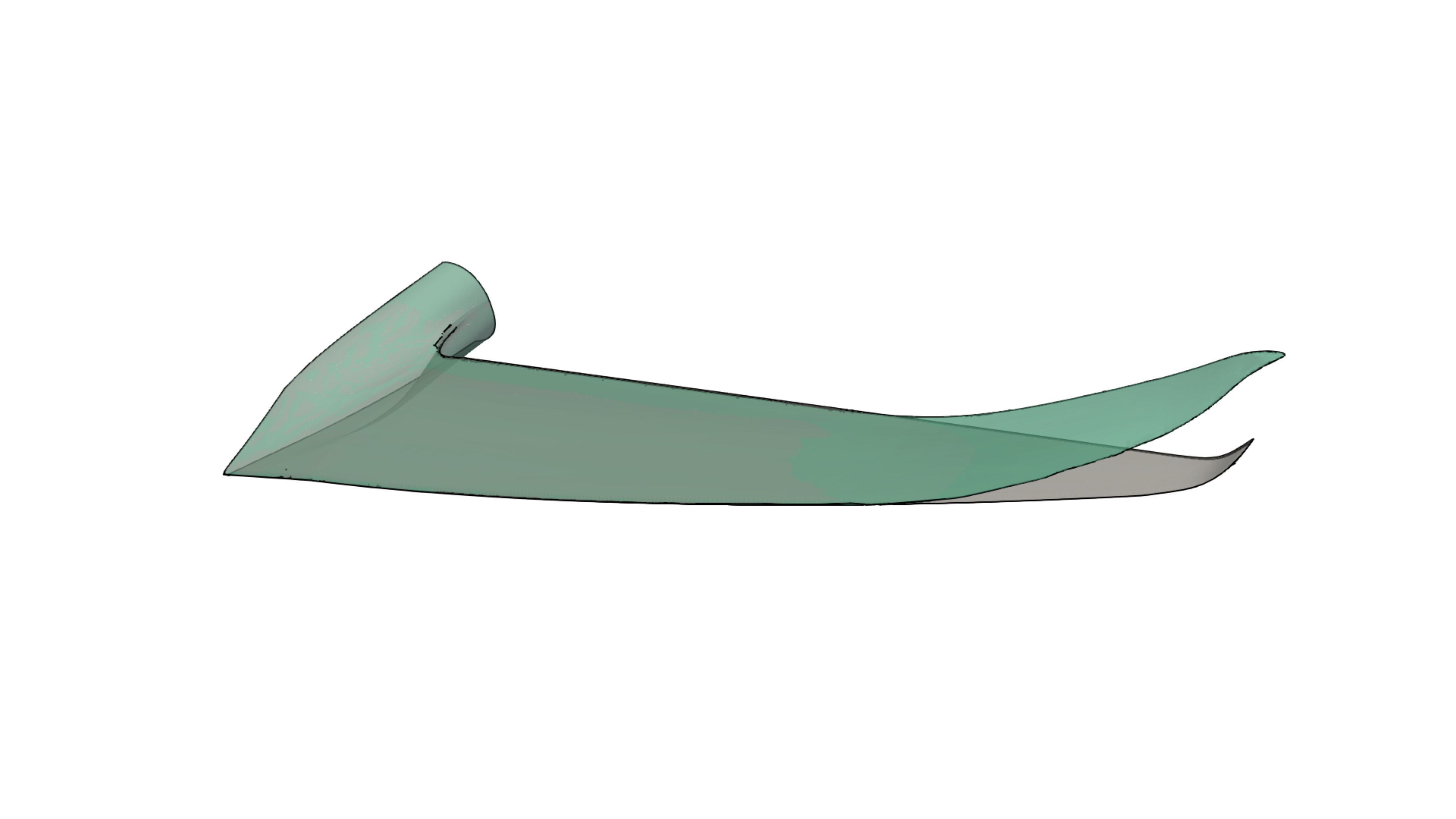}
        \caption{}
    \end{subfigure}
    \hfill
    \begin{subfigure}[b]{0.495\linewidth}
        \centering
        \includegraphics[trim={7.5cm 7.5cm 7.5cm 7.5cm}, clip, width=\linewidth]{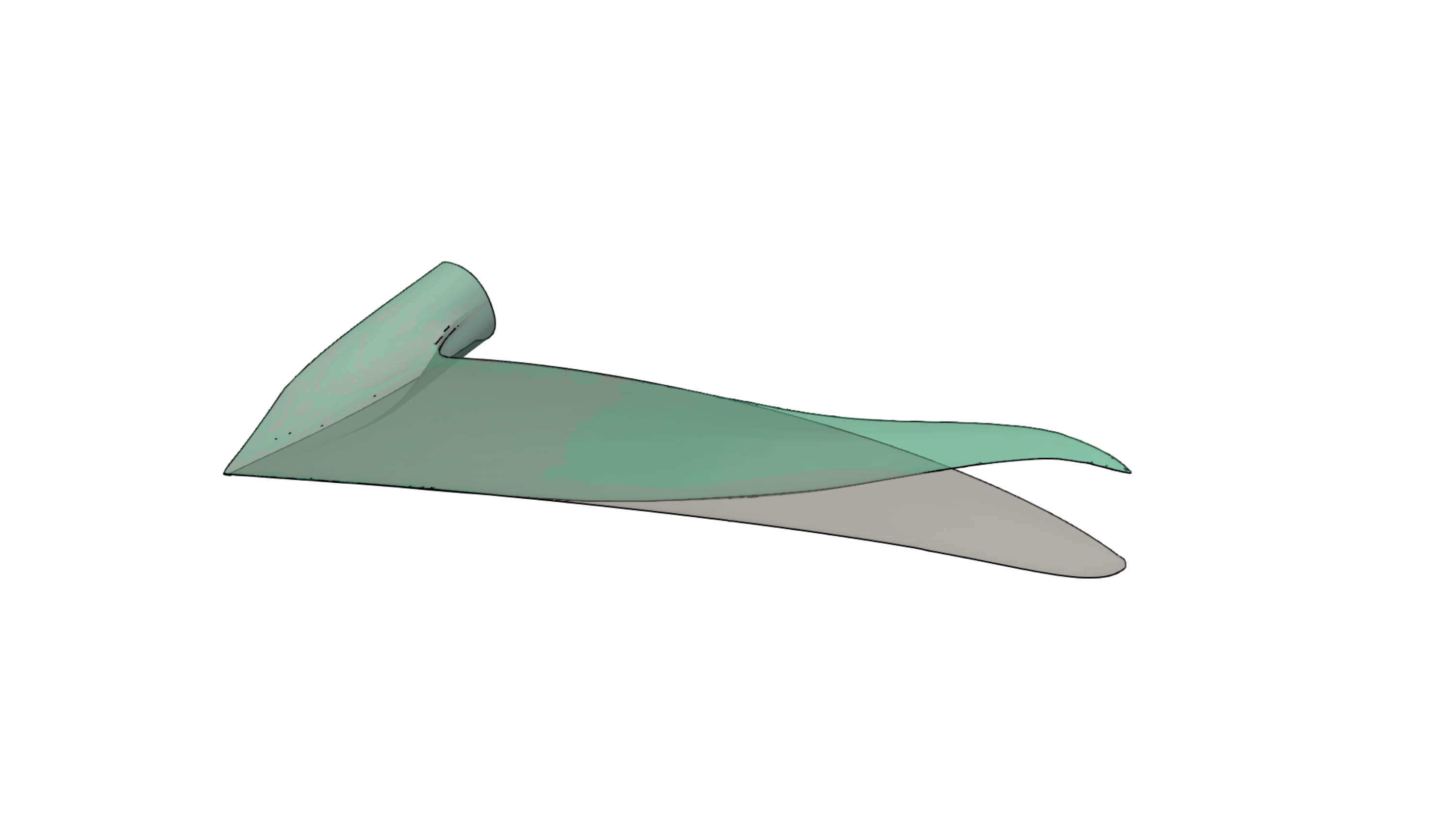}
        \caption{}
    \end{subfigure}
    \caption{Results from the adjoint optimization workflow: (a) best design from the ML-MO workflow versus the adjoint-improved winglet design (green), (b) base geometry versus the adjoint-improved design (green).}
    \label{fig:adjoint_results}
\end{figure}

Figure \ref{fig:adjoint_results} displays geometric modifications achieved through the adjoint optimization process, with the initial designs shown alongside the improved designs (highlighted in green). In both adjoint optimization scenarios, the wing undergoes partial upward deformation, resulting in a shape that differs significantly from the baseline design, which is angled downward. Furthermore, near the tip region, the wing is progressively deflected rearward. These geometric modifications are consistent with the adjoint sensitivity maps, suggesting a design tendency to reduce induced drag by alleviating downwash effects and adjusting the effective loading distribution along the span.

\begin{figure}[H]
	\centering
	\includegraphics[width=\linewidth]{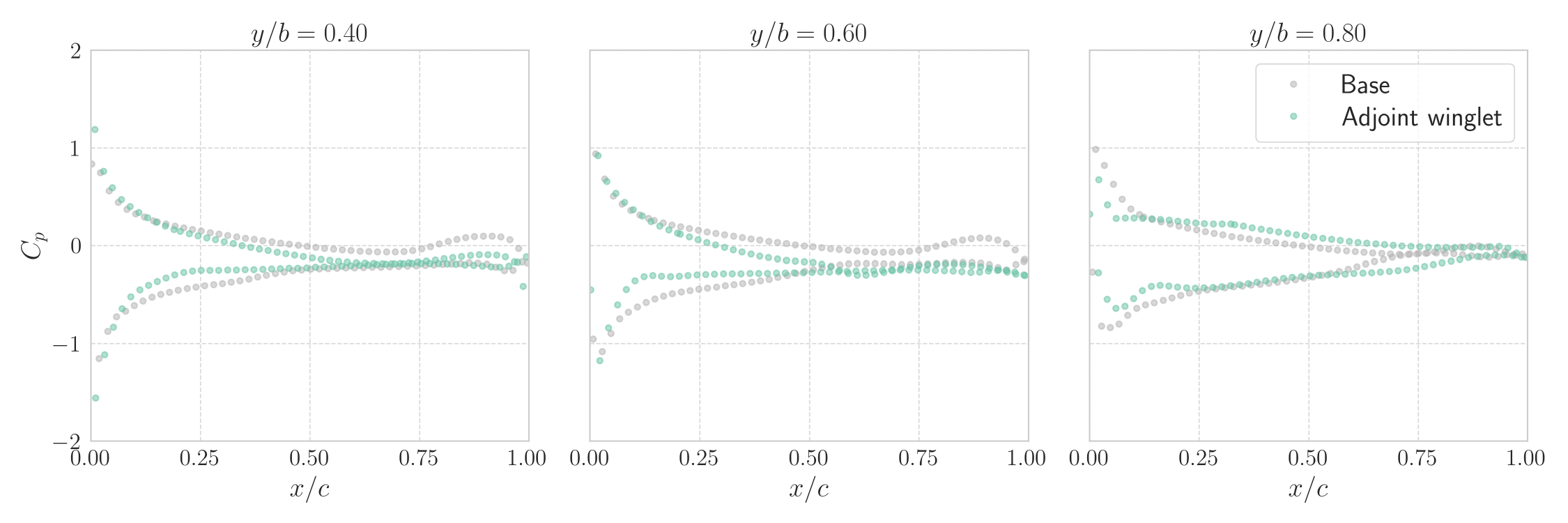}
	\caption{Pressure coefficient $C_p$ at spanwise sections $y/b=0.4$, $0.6$, and $0.8$ (left to right) for the base and adjoint-improved winglet design.}
	\label{fig:adjoint_doe_winglet_results_Cp}
\end{figure}

Pressure coefficient distributions at spanwise locations $y/b = 0.4$, $y/b = 0.6$, and $y/b = 0.8$ are shown in Figure \ref{fig:adjoint_doe_winglet_results_Cp}. The graphs compare the base wing with the adjoint-improved winglet design. The results indicate a deliberate spanwise redistribution of loading in the optimized design. At $y/b = 0.4$, the base design exhibits higher sectional lift. At $y/b = 0.6$, the separation in lift is even larger. Moving outward along the span, at $y/b = 0.8$, the trend reverses, with the adjoint-improved design achieving higher sectional lift.

\begin{figure}[H]
	\centering
	\includegraphics[width=\linewidth]{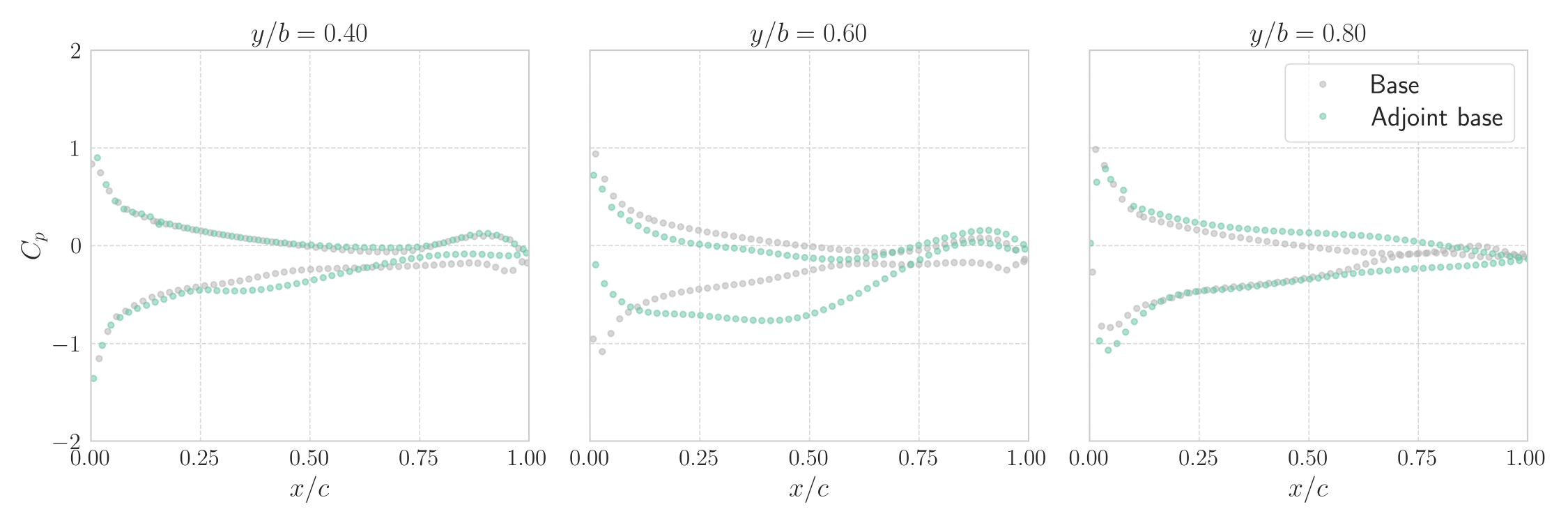}
	\caption{Pressure coefficient $C_p$ at spanwise sections $y/b=0.4$, $0.6$, and $0.8$ (left to right) for the base and adjoint-improved design.}
	\label{fig:adjoint_base_results_Cp}
\end{figure}

Pressure coefficient distributions comparing the base and adjoint-improved base design are shown in Figure \ref{fig:adjoint_base_results_Cp}. The adjoint-improved base design consistently produces higher sectional lift across all spanwise sections. The increase is modest at $y/b = 0.4$ and grows progressively along the span.

\section{Discussion}

The initial analysis revealed that five proposed preliminary configurations exhibited remarkably similar lift-to-drag ratios across different angles of attack. Given that these designs were selected arbitrarily and share similar parameter choices and overall design logic, this performance similarity is not unexpected. The winglet variants aimed to reduce induced drag by modifying tip vortex formation, similar to aircraft wing design principles. However, the arbitrary geometric modifications resulted in suboptimal spanwise load distributions and increased wetted surface area, which offset potential induced benefits. The flipped and split tip designs altered the effective aspect ratio and tip flow characteristics; however, without systematic optimization, these changes did not translate into net performance gains. This also implies that any of the preliminary designs, if properly optimized, could achieve even higher performance. The winglet configuration was optimized to test this hypothesis.

The proposed ML-MO optimization approach, while effective, can be significantly improved. The primary limitation in this study is the number of design parameters, which constrains design flexibility. Expanding the parameter set would substantially increase the search space and could lead to better-performing configurations. However, an exhaustive exploration of an expanded parameter set could become computationally prohibitive. Methods such as Latin hypercube sampling can efficiently cover it with relatively few simulation runs, providing sufficient data for surrogate training. Still, as the number of design variables increases, the optimization effort grows roughly exponentially, making convergence slower and the search for optimal design more challenging.

The adjoint-based approach is a compelling alternative to conventional surrogate-based design methodologies. Unlike parametric optimization that constrains the design space to predefined variables, the adjoint method enables direct modification of the computational domain through gradient-based shape deformation, allowing the optimizer to explore unconventional geometries. This approach offers inherent flexibility by treating the entire boundary as a design variable. However, implementation challenges arise when large geometric deformations occur during the optimization process. The method's effectiveness is also critically dependent on the proper formulation of objective functions and constraints; poorly defined optimization problems can result in infeasible or physically unrealistic solutions. The adjoint approach offers computational advantages over surrogate-based methods. Although a single adjoint simulation took roughly five times longer than a conventional simulation used for DOE, this overhead is compensated by avoiding the extensive simulation campaign that was needed to generate training data. As the number of design parameters for DOE increases, the advantages (efficiency) of the adjoint-based approach become more pronounced.

It must be noted that the study employed several important simplifications, which limit its applicability to real-world conditions. Most notably, all simulations were conducted under incompressible, single-phase, steady-state assumptions. In competitive racing, hydrofoils experience complex multiphysics phenomena. Under extreme hydrodynamic loads, structural deformations can alter the effective geometry and consequently performance characteristics. The influence of unsteady flow effects and wave interactions should also be considered. Incorporating these effects through multiphysics simulations in future work would provide more accurate performance predictions and could reveal design sensitivities not captured in the current idealized approach. Additionally, both optimization methodologies were implemented for a single angle of attack (0$^{\circ}$), which is an idealized case. Optimizing across multiple operating points would improve performance throughout the hydrofoil's typical operating range.

Beyond hydrodynamic and aerodynamic performance, practical factors must also be taken into account. Manufacturing complexity presents a significant constraint, particularly for commercial applications where production costs must remain viable. Complex geometries that offer theoretical performance advantages could prove economically unfeasible or technically challenging to manufacture with sufficient precision and consistency. User experience is another key consideration. Designs optimized for competitive performance may exhibit handling characteristics that are difficult or unsafe for novice sailors, limiting their appeal and market relevance.

\section{Conclusion}

This study presents a systematic methodology for hydrofoil design optimization that integrates reverse engineering, computational fluid dynamics, machine learning surrogate modeling, and adjoint-based optimization. The approach addresses the practical challenge of improving existing hydrofoil designs through computational methods rather than relying solely on empirical adjustments based on intuition and experience.

Two distinct optimization paths were successfully implemented: a data-driven approach using a machine learning surrogate coupled with a multi-objective genetic algorithm, and a gradient-based adjoint optimization approach. The computational efficiency of the surrogate modeling approach makes the methodology accessible for iterative design processes, while the adjoint optimization provides the geometric flexibility needed for discovering novel design concepts. The research revealed that winglet-inspired modifications, while conceptually promising, do not outperform the baseline design across the full operational range.

Several important limitations constrain the generalizability of these results. The study relies on single-phase, steady-state assumptions that do not capture the dynamic behavior characteristic of competitive foil windsurfing. Real-world performance involves complex interactions, wave influences, and structural deformations under dynamic loading conditions. Still, the presented methodology offers significant value for competitive foil windsurfing applications. The 15\% efficiency improvement achieved through adjoint optimization is a substantial achievement that could translate to measurable performance benefits in competitive settings.

\section*{CRediT authorship contribution statement}
\textbf{Ante Sikirica:} Conceptualization, Methodology, Formal analysis, Writing - original draft, Writing - review and editing.
\textbf{Ivana Lučin:} Conceptualization, Methodology, Visualization, Writing - original draft, Writing - review and editing.
\textbf{Marta Alvir:} Methodology, Visualization, Writing - original draft, Writing - review and editing.
\textbf{Tin Matulja:} Writing - original draft, Writing - review and editing, Supervision.

\section*{Declaration of competing interest}
The authors declare that they have no known competing financial interests or personal relationships that could have appeared to influence the work reported in this paper.

\section*{Data availability}
Data will be made available on request.

\section*{Acknowledgements}
This research was partially supported by the UNIRI-INOVA-4-24-7 project of the University of Rijeka. The authors are grateful to the Center for Advanced Computing and Modelling at the University of Rijeka for their support with computational resources.

\bibliographystyle{elsarticle-harv}

\end{document}